\documentclass[10pt]{article}
\usepackage{graphicx}
\usepackage[all,poly]{xy}
\title{The conjugate prior for discrete hierarchical log-linear models }
\author{Jinnan Liu\thanks{Department of Statistics, York
University, Toronto, M3J 1P3, Canada} , H\'el\`ene Massam\thanks{Department of Statistics, York
University, Toronto, M3J 1P3, Canada. This author thanks NSERC
whose generous support  has made this work possible.}}
\date{}
\usepackage{latexsym}
\begin{document}
\maketitle

\def\reel{\hbox{I\hskip-2pt R}}
\def\comp{\hbox{I\hskip-6pt C}}
\def\quat{\hbox{I\hskip-2pt H}}
\def\gene{\hbox{I\hskip-2pt D}}
\def\natu{\hbox{I\hskip-2pt N}}
\def\tr{\;\mathrm{tr}\;}
\def\<{\langle}
\def\>{\rangle}
\def\P{{\bf P}}
\def\E{{\bf E}}
\def\n{{\boldmath n}}
\def\b{{\boldmath b}}
\def \m{{\bf m}}
\def\z{{\boldmath z}}

\def\<{\langle}
\def\>{\rangle}

\def\th{th\'eor\`eme\ }
\def\ths{th\'eor\`emes\ }
\def\reel{\hbox{I\hskip -2pt R}}
\def\pe{\hbox{I\hskip -2pt P}}
\def\u{{\bf u}}
\def\v{{\bf v}}
\def\esp{{\bf E}}
\def\<{\langle}
\def\>{\rangle}
\def\x{{\bf x}}
\def\y{{\bf y}}
\def\f{{\bf f}}
\def\p{{\bf p}}
\def\w{{\bf w}}
\def\s{{\bf s}}
\def\W{{\bf W}}
\def\Z{{\bf Z}}
\def\Q{{\bf Q}}
\def\tr{\textmd{trace}\,}
\def\grad{\textmd{grad}\,}

\textheight =22.5 cm
\textwidth =15 cm
\voffset =-0.5 in
\hoffset =0 in
\headheight =0 cm
\renewcommand{\thesection}{\arabic{section}}
\renewcommand{\theequation}{\arabic{section}\mathrm{.}\arabic{equation}}

\newtheorem{theorem}{Theorem}[section]
\newtheorem{definition}{Definition}[section]
\newtheorem{subsecdef}{Definition}[subsection]
\newtheorem{lemma}{Lemma}[section]
\newtheorem{prop}{Proposition}[section]
\newtheorem{remarks}{Remarks}[section]
\newtheorem{remark}{Remark}[section]
\newtheorem{example}{Example}[section]
\newtheorem{perty} {Property}[section]
\newtheorem{cor} {Corollary}[section]

\newenvironment{pff}{\hspace*{-\parindent}{\bf Proof:}}{\hfill $\Box$
\vspace*{0.2cm}}
\newenvironment{pff.}{\hspace*{-\parindent}{\bf Proof}}{\hfill $\Box$
\vspace*{0.2cm}}
\renewcommand\theequation{\thesection .\arabic{equation}}
\let\subs\subsection

\let\sect\section
\renewcommand\section[1]{\setcounter{equation}{0}\sect{#1}}

\newcommand{\pperp}{{\ \perp\!\!\!\perp\ }}
\def\tr{\;\mathrm{tr}\;}
\def\Si{\Sigma^{\cal V}}

\newcommand{\cK}{{\cal K}}
\newcommand{\cM}{{\cal M}}
\newcommand{\IR}{{\Bbb R}}
\newcommand{\IC}{{\Bbb C}}
\newcommand{\ID}{{\Bbb D}}
\newcommand{\IO}{{\Bbb O}}
\newcommand{\IN}{{\Bbb N}}
\newcommand{\IH}{{\Bbb H}}
\begin{abstract}
In the Bayesian analysis of contingency table data, the selection of a
prior distribution for either the log-linear parameters or the cell
probabilities parameter is a major challenge.
Though the conjugate prior on cell probabilities
 has been defined by Dawid and Lauritzen (1993) for decomposable graphical models, it 
has not been identified for the  larger class of graphical models Markov with respect to an arbitrary undirected graph
or for the even wider class of hierarchical log-linear models. In this paper, working with the log-linear parameters used by GLIM, we first define the conjugate prior for  these parameters and then derive the induced prior for the cell probabilities: this is done for the general class of hierarchical log-linear models.
We show that the conjugate  prior has all the properties that one expects from a prior: notational simplicity, ability to reflect either no prior knowledge or a priori expert knowledge, a moderate number of hyperparameters and mathematical convenience. It also has the strong hyper Markov property which allows for local updates within prime components for graphical models.

\vspace{.5cm}

Keywords: hierarchical log-linear models, conjugate prior, prior specification, hyper Markov property.
\vspace{.5cm}

AMS 2000 Subject classifications. Primary 62H99; Secondary 62E15
\end{abstract}

\newpage
\section{Introduction}

We consider data given under the form of a contingency table representing the classification of $n$ individuals according to a finite set $V$ of criteria. Each criterion $\gamma\in V$ is represented by a variable $X_{\gamma}$ which take values in a finite set ${\cal I}_{\gamma}$. We observe the values of the variable $X=(X_{\gamma}, {\gamma}\in V)$ in ${\cal I}=\times_{{\gamma}\in V}{\cal I}_{\gamma}$ for these $n$ individuals,  and we assume that the resulting $|{\cal I}|$ cell counts in the contingency table follow a multinomial distribution. We also assume that the cell probabilities are modeled according to a hierarchical log-linear model parametrized by interaction parameters represented by a vector $\theta$. Since the class of discrete graphical models Markov with respect to an arbitrary undirected graph $G$ is an important subclass of the class of hierarchical log-linear models, we will give special attention to that class througout the paper.

In the Bayesian analysis of contingency table data, the selection of a
prior distribution for either the log-linear parameters or the cell
probabilities parameter is a major challenge (see Clyde and George, 2004). Priors are usually chosen for their conceptual and computational simplicity and for their ability to represent experts prior beliefs. They are also chosen so that they can conveniently be used for  the whole class of log-linear models which includes nondecomposable as well as decomposable graphical models. Moreover their parametrization, that is the hyper-parametrization, should be such that
hyper-parameters are compatible across models.

As shown in Dawid and Lauritzen (1993), the conjugate prior for decomposable graphical models, called the hyper Dirichlet and defined for marginal cliques and separators  cell probabilities has all of these properties and additionally has the strong hyper Markov property. The latter is very desirable since it allows for local updates within cliques thus simplifying the computation of Bayes factors in a model selection process. The hyper Dirichlet has therefore been used in many studies( see for example Madigan and Raftery, 1994 and Dellaportas and Forster, 1999).
However, the hyper Dirichlet is only defined for decomposable graphical models and when it is used as a prior, the corresponding posterior probability for a model is only its probability within this restricted class thus making it difficult to compare it to the posterior probability of another model considered within the wider class of hierarchical log-linear models. Moreover, it appears to have many hyperparameters since a set of parameters has to be chosen for the Dirichlet on each clique and each separator of the graph. In fact, all these hyper-parameters are not independent of each other since they have to be hyper-consistent but the apparently large number of parameters adds a level of complexity to their selection.

 Consequently much effort has been devoted to the study of alternative priors. For  example, King and Brooks (2001), after a discussion on the advantages and disadvantages of the hyper Dirichlet propose a multivariate normal prior for the log-linear parameters for all hierarchical log-linear models. This prior allows for efficient computation, facilitates prior elicitation and induces a log-normal distribution on the cell probabilities with easy to compute prior mean and covariances.

The aim of this paper is to show that  the conjugate prior can also be defined for the wider class of hierarchical log-linear models in a simple way and  that it has all the desirable properties that one traditionally wants from a prior. Indeed we will show that  experts prior beliefs or lack of any prior information can easily be expressed by an appropriate choice of hyperparameters. The chosen prior is consistent with prior beliefs under both parametrization of the model.
The conjugate prior is also hyper Markov thus leading to local updates in graphical models, a property that traditional normal priors on log-linear parameters do not have. Also, the number of hyperparameters is moderate, in fact exactly equal to the number of log-linear parameters plus one and the hyperparameters are hyperconsistent across prime components in graphical models and compatible across models.

In \S2, we set our notation and give some preliminary results. We work with the parametrization used by GLIM. It is interesting to note that this parametrization expresses the logarithm of cell probabilities $p(i), i\in {\cal I}$, which we regard as functions of $i$, as the sum of functions $\theta_E(i)$ of $i$ which are in orthogonal subspaces of the space $\reel^{\cal I}$  of functions on ${\cal I}$.
This orthogonal decomposition of $\log p(i)$ will insure that hyperparameters in the prior are compatible across all models.
Our parametrization will also lead us to express the distribution of the marginal cell counts in the contingency table, rather than the cells counts, as an exponential family. Using this exponential family form, we derive, in \S3, the expression of the  Diaconis and Ylvisaker (1979) conjugate prior for the log-linear parameters. We give a necessary and sufficient condition for this prior to be proper and two methods to obtain hyperparameters that insure that the prior is proper. In \S4, we obtain the expression of the induced conjugate prior for the cell probabilities and in \S5, we give the details of the properties we mentioned above. Having the expression of the induced prior on cell probabilities allow us to verify that 
the choice of hyperparameters in one parametrization (log-linear or cell probabilities) expresses the same prior belief in the other parametrization.

\section{The log-linear model}
\subsection{The parametrization}
Let $V$ be the set of criteria. Let $X=(X_{\gamma},|\;\gamma\in V)$ such that $X_{\gamma}$ takes its values (or levels) in the finite set $I_{\gamma}$ of dimension $|I_{\gamma}|$. When a fixed number of individuals are classified according to the $|V|$ criteria, the data is collected in a contingency table with cells indexed by combination of levels for the $|V|$ variables. We adopt the notation of Lauritzen (1996) and denote a cell by
$$i=(i_{\gamma},\;{\gamma}\in V)\in {\cal I}=\times_{{\gamma}\in V}{\cal I}_{{\gamma}}.$$
The count in cell $i$ is denoted $n(i)$ and the probability of an individual falling in cell $i$ is denoted $p(i)$.
For $E\subset V$, cells  in the $E$-marginal table are denoted
$$i_E\in {\cal I}_E=\times_{{\gamma}\in E}{\cal I}_{{\gamma}}.$$
The marginal counts are denoted $n(i_E).$
For $n=\sum_{i\in {\cal I}}n(i),\;(n)=(n(i),\;i\in {\cal I})$ follows a multinomial ${\cal M}(n, p(i),\;i\in {\cal I})$ distribution with probability density function
\begin{equation}
\label{pofn}
P(n)\propto \prod_{i\in {\cal I}}p(i)^{n(i)}\;.
\end{equation}
Let $i^*$ be a fixed but arbitrary cell which for convenience we take to be the cell indexed by the "lowest levels" for each factor and for convenience again, we denote this level by 0. Therefore $i^*$ can be thought to be the cell $$i^*=(0,0,\ldots,0).$$

Consider the following parametrization 
\begin{equation}
\label{pitotheta}
\theta_E(i)=\sum_{F\subseteq E}(-1)^{|E\setminus F|}\log p(i_F,i^*_{F^c})
\end{equation}
where by Moebius inversion
\begin{equation}
\label{probi}
p(i)=\exp \sum_{E\subseteq V}\theta_E(i)\;.
\end{equation}
We note that $\theta_{\emptyset}(i)=\log p(i^*), \;i\in {\cal I}$ and we will therefore adopt the notation
\begin{equation}
\label{pto}
\theta_{\emptyset}(i)=\theta_{\emptyset},\;\;p(i^*)=p_{\emptyset}=\exp \theta_{\emptyset}.
\end{equation}
This parametrization  has been used in many papers (see for example Dellaportas and Forster, 1999) and can be found in Lauritzen (1996, p.36).

 Let us make an important remark here. Since $i^*$ is fixed, the function
$$i\in {\cal I}\mapsto \log p(i_E, i^*{E^c})$$
belongs to the factor subspace ${\cal U}_E$ ( as defined in Lauritzen ,1996, Appendix B.2) of the space $\reel^{\cal I}$ of real-valued function on ${\cal I}$ that depend only on $i_E$. Therefore by Proposition B.4 of Lauritzen (1996) and 
(\ref{pitotheta}) above, $\theta_E(i)$ belongs to the interaction subspace ${\cal V}_E$ which gives the "`pure"' contribution of the interaction between variables in $E$ with the interaction between variables in all $F\subset E$ removed. This means that (\ref{probi}) or more precisely its equivalent expression
$$\log p(i)=\sum_{E\subseteq V}\theta_E(i)$$ 
is the unique expansion of $\log p(i)$ into its orthogonal components in ${\cal V}_E, E\subseteq V$. This orthogonal decomposition of $\log p(i)$ is the property that will make the hyperparametrization of the conjugate prior on $\theta$, defined below in (\ref{freethetadiscrete}), 
compatible across models since all models will be expressed in the same orthogonal "basis".
Let us now emphasize some other properties of the $\theta$ parametrization with the following three lemmas.

\begin{lemma}
\label{thetaona}
For any $(i)\in {\cal I}$ and any $E\subseteq V$, $\theta_E(i)$ depends only on $i_E$, that is
$$\theta_E(i)=\theta_E(i_E)\;.$$
\end{lemma}
Since  $i^*$ is fixed, the proof of this first lemma is obvious.

\begin{lemma}
\label{glim}
If $i$ is such that for $\gamma \in E, i_{\gamma}=i^*_{\gamma}=0$, then $\theta_E(i_E)=0$
\end{lemma}
\begin{pff}
By definition and since $(i_{F\cup \gamma},i^*_{(F\cup \gamma)^c})=(i_{F},i^*_{F^c})$ we have
\begin{eqnarray}
\theta_E(i)&=&\sum_{F\subseteq E\setminus \gamma}(-1)^{|E\setminus F|}\log p(i_F,i^*_{F^c})-\sum_{F\subseteq E\setminus \gamma}(-1)^{|E\setminus F|}\log p(i_{F\cup \gamma},i^*_{(F\cup \gamma)^c})\nonumber\nonumber\\
&=&\sum_{F\subseteq E\setminus \gamma}(-1)^{|E\setminus F|}\log p(i_F,i^*_{F^c})-\sum_{F\subseteq E\setminus \gamma}(-1)^{|E\setminus F|}\log p(i_{F},i^*_{F^c})=0\;.\nonumber
\end{eqnarray}
\end{pff}

\noindent From this lemma, it follows immediately that our parametrization is the GLIM parametrization that sets to 0 the values of the $E-$ interaction log-linear parameters when at least one index in $E$ is at level  0 (see for example Agresti 1990, p.150) . Therefore, for each $E\subseteq V$, there are only $\prod_{\gamma \in E}(|{\cal I}_{\gamma}|-1)$ parameters. 
 The next lemma is actually the Hammersley-Clifford theorem and its proof can be found for example in Lauritzen (1996, p.36).
\begin{lemma}
\label{hamcli}
Assuming all cell probabilities are positive, the distribution of $X=(X_{\gamma}, \gamma\in V)$ is Markov with respect to the undirected graph $G$ if and only if $\theta_E(i_E)=0$ whenever  $E\subseteq V$ is not complete. 
\end{lemma}

From this lemma, it follows that the multinomial distribution of the cell counts is Markov with respect to a graph $G$ if and only if $\theta_E$ is equal to zero when $E$ is not complete, a well-known property that we recall here (see Darroch, Lauritzen and Speed, 1980).

 \noindent For notational convenience, we now define
$${\cal E}=\{E\subseteq V,\;\; E\not = \emptyset\}.$$
By Lemma \ref{glim}, for any given $j\in {\cal I}$, $\theta_E(j_E)=0$ if there exists at least one $\gamma\in E$ such that  $j_{\gamma}=0$. We therefore define for any $j\in {\cal I}$
\begin{equation}
{\cal E}^*_j=\{E\in {\cal E}\;|\;j_{\gamma}\not = 0, \forall \gamma\in E   \}\;.
\end{equation}
Then, by (\ref{pitotheta}), (\ref{pto}) and Lemma (\ref{thetaona}), 
\begin{eqnarray}
p_{\emptyset}&=&1-\sum_{i\in {\cal I}, i\not = i^*}p(i)=1-\sum_{i\in {\cal I}, i\not = i^*}\exp  \sum_{E\subseteq V}\theta_E(i)
=1-\sum_{i \in {\cal I}, i\not = i^*}p_{\emptyset}\Big(\exp\sum_{E\in {\cal E}^*_i}\theta_E(i_E)\Big)\nonumber\;.
\end{eqnarray}
which yields 
\begin{equation}
\label{piempty}
p_{\emptyset}=\frac{1}{1+\sum_{j \in {\cal I}, j\not = i^*}\exp\sum_{E\in {\cal E}^*_j}\theta_E(j_E)}\;\;
\end{equation}
and
\begin{equation}
\label{pigeneral}
p(i)=\frac{\exp \sum_{E\in {\cal E}^*_i}\theta_E(i_E)}{1+\sum_{j \in {\cal I}, j\not = i^*}\exp\sum_{E\in {\cal E}^*_j}\theta_E(j_E)}
\end{equation}
and thus all cell probabilities are expressed in terms of the free parameters $$\theta_E(i_E), E\in {\cal E}, i\in {\cal I}, i\not = i^*.$$
\subsection{The multinomial distribution for discrete data}
We now want to give the probability density function of the multinomial distribution under the form of an exponential family. This will be done successively for the saturated model i.e. Markov with respect to a complete graph, for models Markov with respect to an undirected graph $G$ and for general hierarchical log-linear models.

\noindent From  (\ref{pofn}), for the saturated model, we have
\begin{eqnarray}
P(n)&\propto& \prod_{i\in {\cal I}}(\exp \sum_{E\subseteq V} \theta_E(i))^{n(i)}=\prod_{i\in {\cal I}}\exp n(i)\log \exp \sum_{E\subseteq V} \theta_E(i_E)=\exp \sum_{i\in {\cal I}} n(i)\sum_{E\subseteq V} \theta_E(i)\nonumber\\
&=&\exp \sum_{i\in {\cal I}}\sum_{E\subseteq V}\theta_E(i_E)\sum_{j\in {\cal I}, j_E=i_E} n(j)=\exp \sum_{E\subseteq V}\sum_{i_E\in {\cal I}_E}\theta_E(i_E)n(i_E)\nonumber\\
&=&\exp \{\sum_{E\subseteq V, E\not = \emptyset}\sum_{i_E\in {\cal I}_E}\theta_E(i_E)n(i_E)+n\theta_{\emptyset)}\}\nonumber
\end{eqnarray}
Moreover, we know from Lemma \ref{glim} that, only those $\theta_E(i_E)$ where $i_{\gamma}\not =0, \gamma\in E$ are nonzero. Therefore if, for $E\in {\cal E}$ we define
\begin{equation}
{\cal I}_E^*=\{i_E=(i_{\gamma}, \gamma\in E)\in {\cal I}_E,\; i_{\gamma}\not =0, \gamma\in E    \}
\end{equation}
then the probability density function above becomes
\begin{equation}
P(n)\propto \exp \{\sum_{E\subseteq V, E\not = \emptyset}\sum_{i_E\in {\cal I}_E^*}\theta_E(i_E)n(i_E)+n\theta_{\emptyset)}\}
\end{equation}

We see that, with the parametrization that we have chosen, the marginal counts $n(i_E)$, rather than the cell counts $n(i)$, appear naturally as random variables .
Since the Jacobian of
$$n=(n(i), i\in {\cal I})\mapsto y=(n(i_E),\; E\in {\cal E},\; i_E\in {\cal I}_E^*)$$
is clearly one,
 the family of  distributions for $y$ is the natural exponential family
\begin{equation}
\label{nefi}
{\cal F}_{\mu}=\{f(y; \theta)\mu(y)= \frac{\exp \{\sum_{E\in {\cal E}}\sum_{i_E\in {\cal I}_E^*}\theta_E(i_E)n(i_E)\}
}{\Big(1+\sum_{i \in {\cal I}, i\not = i^*}\exp \sum_{F\in {\cal E}^*_i}\theta_F(i_F)\Big)^n}\mu(y) ,\;\; \theta \in  \reel^{\prod_{E\in {\cal E}}\prod_{\gamma\in E}(|I_{\gamma}|-1)}  \}\;,
\end{equation}
where $\mu$ is a reference measure of no particular interest to us here. This gives us the density for the saturated model.
\vspace{2mm}

\noindent When $G$ is an arbitrary undirected graph let
$${\cal D}=\{E\in {\cal E}\;|\; E\;\;\mbox{complete}\}\;$$
and for any given $j\in {\cal I}$
\begin{equation}
{\cal D}^*_j=\{D\in {\cal D}\;|\;j_{\gamma}\not = 0, \forall \gamma\in D   \},
\end{equation}

From Lemma \ref{hamcli}, to obtain the expression of the cell probabilities and of the family of multinomial distributions Markov with respect to $G$ it suffices to equate $\theta_E(i)$ to 0 for $E\not \in {\cal D}$ and for all $(i)$. We therefore have
\begin{equation}
\label{piemptyg}
p_{\emptyset}=\frac{1}{1+\sum_{j \in {\cal I}, j\not = i^*}\exp\sum_{D\in {\cal D}^*_j}\theta_D(j_D)}\;\;
\end{equation}
\begin{equation}
\label{pigeneralg}
p(i)=\frac{\exp \sum_{D\in {\cal D}^*_i}\theta_D(i_D)}{1+\sum_{j \in {\cal I}, j\not = i^*}\exp\sum_{D\in {\cal D}^*_j}\theta_D(j_D)}\;,
\end{equation}
where it is important to note that in  (\ref{pigeneralg})
not all $p(i)$ are free parameters since 
$$\mbox{for}\;\;\; E\not \in {\cal D},\;\;\; \theta_E(i_E)=0$$
which implies that for $E\not \in {\cal D},\; p(i_E, i^*_{E^c})$ is function of $p(i_F, i^*_{F^c}), F\subset E,\; F\in {\cal D}\;.$
Only cell probabilities of the form
$p(i_D,i^*_{D^c}),\;i_D\in {\cal I}^*_D,\;D\in {\cal D}$
will be free probabilities and form the cell probability parameter
\begin{equation}
\label{freepdiscrete}
p=(p(i_D,i^*_{D^c}),\;D\in {\cal D},\;i_D\in {\cal I}^*_D,)\;\;\mbox{with}\;\;p(i_D, i^*_{D^c})\;\;\mbox{as in}\;\;(\ref{probi})\;.
\end{equation}
of the multinomial distribution Markov with respect to $G$ for graphical models, or of the hierarchical log-linear model.
The corresponding log-linear parameters are obviously
\begin{equation}
\label{freethetadiscrete}
\theta=(\theta_D(i_D),\;D\in {\cal D},\;i_D\in {\cal I}^*_D)\;\;\mbox{with}\;\;\theta_D(i_D)\;\;\mbox{as in}\;\;(\ref{pitotheta})\;.
\end{equation}
\noindent Moreover, the family of multinomial distribution Markov with respect to $G$ for $$y=(n(i_D),\;D\in {\cal D},\; i_D\in {\cal I}_D^*)$$ is 
\begin{eqnarray}
{\cal F}_{\mu_G}&=&\{f_G(y; \theta)\mu_G(y)= \frac{\exp \{\sum_{D\in {\cal D}}\sum_{i_D\in {\cal I}_D^*}\theta_D(i_D)n(i_D)\}
}{\Big(1+\sum_{i \in {\cal I}, i\not = i^*}\exp \sum_{D\in {\cal D}^*_i}\theta_D(i_D)\Big)^n}\mu_G(y) ,\nonumber\\
\label{nefig}
&&\hspace{8cm} \theta \in \reel^{\prod_{D\in {\cal D}}\prod_{\gamma\in D}(|{\cal I}_{\gamma}|-1)}  \}\;,
\end{eqnarray}
where $\theta=(\theta_D(i_D), D\in {\cal D}, i_D\in {\cal I}_D)$ and 
$\mu_G$ is a reference measure of no particular interest to us here. Densities in ${\cal F}_{\mu_G}$ will be written under the natural exponential family form
\begin{eqnarray}
\label{fg}
f_G(y; \theta)&=&\exp \{\sum_{D\in {\cal D}}\sum_{i_D\in {\cal I}_D^*}\theta_D(i_D)n(i_D)\}-n\log \Big(1+\sum_{i \in {\cal I}, i\not = i^*}\exp \sum_{D\in {\cal D}^*_i}\theta_D(i_D)\Big)
\end{eqnarray}

\vspace{2mm}

When the model is a hierarchical log-linear model, let ${\cal D}$ be the set of subsets of $V$ representing  the set of
all possible interactions in the given model, which we will call the generating set. Then, the expression of the cell probabilities and of the multinomial distribution for this model is the same as in (\ref{piemptyg}), (\ref{pigeneralg}) and (\ref{nefig}) but with ${\cal D}$ representing the generating set for the model.

\subsection{The multinomial distribution for binary data}
We consider here the important special case of binary data, because it occurs often in practice and also because in this case, the notation is somewhat simpler.
When the variables $X_{\gamma}, \gamma\in V$ can only take two values 0 or 1, there is only one cell $i$ in each ${\cal I}_E^*, E\in {\cal E}$  and therefore each cell $i=(i_{\gamma}, \gamma\in V)$ can be indexed by 
$E=\{\gamma\in V: i_{\gamma}=1\}$ for  $E\in {\cal E}\cup \emptyset$ such that $i_{\gamma}\not = 0, \gamma\in E$.
The correspondence between ${\cal I}$ and ${\cal E}\cup \emptyset$ is one to one. 
For $i=(i_{\gamma}=1, \gamma\in E, i_{\gamma}=0, i\not \in E)$, we will therefore use the notation 
$$p_E=p(i)\;\;\mbox{for}\;\;\mbox{and}\;\;\theta_E=\theta _E(i_E).$$

The relation (\ref{pitotheta}) becomes
\begin{equation}
\label{thetae}
 \theta_F=\sum_{E\subseteq F}(-1)^{|F\setminus E|}\log p_E=\log \prod_{E\subseteq F}p_E^{(-1)^{|F\setminus E|}},\;\;\; F\in {\cal E}\cup \emptyset
 \end{equation}
\noindent Or equivalently by Moebius inversion

\begin{equation}
\label{ptotheta}
\log p_F=\sum_{E\subseteq F}\theta_E,\;\; F\in {\cal E}\;\;\;\mbox{with} \;\;\;\log p_{\emptyset}=\theta_{\emptyset}\;,
\end{equation}

Then (\ref{piempty}), (\ref{pigeneral}) and (\ref{nefig}) become respectively
\begin{eqnarray}
p_{\emptyset}&=&\frac{1}{\Big(1+\sum_{E\in {\cal E}}\exp \{ \sum_{F\subseteq E, F\in {\cal D}}\theta_F\}\Big)}\label{newp0}\\
p_E&=&\frac{
\exp \sum_{F\subseteq E, F\in {\cal D}}\theta_F}
{\Big(1+\sum_{H\in {\cal E}}\exp \{ \sum_{F\subseteq H, F\in {\cal D}}\theta_F\}\Big)},\;\;E\in {\cal E}
\label{pe}
\end{eqnarray}
and
\begin{eqnarray}
\label{nefmarkov}
{\cal F}_{\mu_G}&=&\{
f(y;\theta, G)\mu_G(y)=\exp\;\Big(
\sum_{D\in {\cal D}}\theta_Dy_D-n\log (1+\sum_{E\in {\cal E}}\exp (\sum_{D\subseteq E, D\in {\cal D}}\theta_D))\Big)\;\mu_G(y)\\
&&\hspace{9cm}\theta=(\theta_D, D\in {\cal D})\in \Theta_G=\reel^{|{\cal D}|}\}\nonumber
\end{eqnarray}
where ${\cal D}$ is equal to ${\cal E}$, the set of complete subsets of $V$ in $G$ or the generating set for the hierarchical model for, respectively, the saturated model, graphical model with respect to $G$ or the hierarchical model.

\noindent 
We note that (\ref{freepdiscrete}) and (\ref{freethetadiscrete}) become
$$p=(p_D, D\in {\cal D})\;\;\mbox{and}\;\;\theta=(\theta_D, D\in {\cal D})$$
respectively.

\subsection{An example}
\noindent  We consider the case where $X=(X_a,X_b,X_c,X_d)$ is Markov with respect to the four-cycle as given below and where the variables are binary.

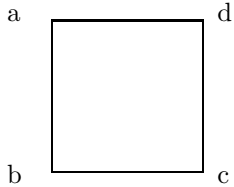
\begin{figure}[h]
\vspace{1cm} \setlength{\unitlength}{1cm}
\begin{picture}(2,2)
\put(5.40,2){a} \put(8.2,2){d} \put(5.40,-0.15){b}
\put(8.2,-0.15){c} \put(6,0){\framebox(2,2){}}
\label{firstfigure}
\end{picture}
\end{figure}
\vspace{4mm}

We then have
\begin{eqnarray}
{\cal D}&=&\{a,b,c,d,ab,bc,cd,da\}\nonumber\\
{\cal E}&=&\{a,b,c,d,ab,bc,cd,da,ac,bd,abc,bcd,cda,dab,abcd\}\nonumber
\end{eqnarray}

The linear constraints on $\theta_E, E\not \in {\cal D}$ are
\begin{eqnarray}
&&\theta_{ac}=\theta_{bd}=\theta_{abc}=\theta_{bcd}=\theta_{cda}=\theta_{dab}=\theta_{abcd}=0\nonumber
\end{eqnarray}
and using (\ref{thetae}), we obtain $p_E, E\not \in {\cal D}$ in terms of $p=(p_E,\; E\in {\cal D})$ as follows
\begin{eqnarray}
&&p_{ac}=\frac{p_ap_c}{p_{\emptyset}},\;p_{bd}=\frac{p_bp_d}{p_{\emptyset}},\;
p_{abc}=\frac{p_{ab}p_{bc}}{p_{b}},\;p_{bcd}=\frac{p_{bc}p_{cd}}{p_{c}},\;p_{cda}=\frac{p_{cd}p_{da}}{p_{d}},\;p_{dab}=\frac{p_{da}p_{ab}}{p_{a}}\;,\nonumber\\
&&p_{abcd}=\frac{p_{ab}p_{bc}p_{cd}p_{da}p_{\emptyset}}{p_{a}p_{b}p_{c}p_{d}}\;.\nonumber
\end{eqnarray}
The cell probability parameters of the multinomial distribution Markov with respect to the four-cycle above can be written in terms of $\theta$  as
\begin{eqnarray}
p_{\emptyset}^{-1}&=&1+e^{\theta_a}+e^{\theta_b}+e^{\theta_c}+e^{\theta_d}+e^{\theta_a+\theta_b+\theta_{ab}}+e^{\theta_b+\theta_c+\theta_{ab}}+e^{\theta_c+\theta_d+\theta_{cd}}+e^{\theta_d+\theta_a+\theta_{da}}\nonumber\\
&&\hspace{2cm}+e^{\theta_a+\theta_b+\theta_c+\theta_{ab}+\theta_{bc}}+e^{\theta_b+\theta_c+\theta_d+\theta_{bc}+\theta_{cd}}+e^{\theta_c+\theta_d+\theta_a+\theta_{cd}+\theta_{da}}+e^{\theta_d+\theta_a+\theta_b+\theta_{da}+\theta_{ab}}\nonumber\\
&&\hspace{2cm}+e^{\theta_a+\theta_b+\theta_c+\theta_d+\theta_{ab}+\theta_{bc}+\theta_{cd}+\theta_{da}}\nonumber\\
p_D&=&p_{\emptyset}e^{\theta_D},\;D\in \{a,b,c,d,\}\;\;\;\mbox{with}\nonumber\\
p_{ab}&=&p_{\emptyset}e^{\theta_a+\theta_b+\theta_{ab}},\;
p_{bc}=p_{\emptyset}e^{\theta_b+\theta_c+\theta_{bc}},\;
p_{cd}=p_{\emptyset}e^{\theta_c+\theta_d+\theta_{cd}},\;
p_{da}=p_{\emptyset}e^{\theta_d+\theta_a+\theta_{da}},\nonumber\\
p_{abc}&=&\frac{p_{ab}p_{bc}}{p_b},\;p_{bcd}=\frac{p_{bc}p_{cd}}{p_c},\;p_{cda}=\frac{p_{cd}p_{da}}{p_d},\;p_{dab}=\frac{p_{da}p_{ab}}{p_a},\;
p_{abcd}=p_{\emptyset}\frac{p_{ab}p_{bc}p_{cd}p_{da}}{p_{a}p_{b}p_{c}p_{d}}\nonumber
\end{eqnarray}

\section{The conjugate prior for the log-linear parameter $\theta$}
From (\ref{fg}), it is clear that, for the three nested classes of models considered in this paper, graphical with respect to $G$ decomposable, graphical with respect to an arbitrary undirected $G$  and hierarchical, the probability density function for the marginal counts $y$ can be written under an exponential family form and therefore the form of   the conjugate prior for $\theta$  is given immediately (see Diaconis and Ylvisaker, 1979) by
\begin{eqnarray}
\label{pitheta}
\hspace{-5mm} \pi_G(\theta|s,\alpha)&=&I_G(s,\alpha)^{-1}\exp \{\sum_{D\in {\cal D}}\sum_{i_D\in {\cal I}_D^*}\theta_D(i_D)s(i_D)\}-\alpha\log \Big(1+\sum_{i \in {\cal I}, i\not = i^*}\exp \sum_{D\in {\cal D}^*_i}\theta_D(i_D)\Big)
\end{eqnarray}
where $I_G(s,\alpha)$ is the normalising constant
\begin{equation}
\label{i}
I_G(s,\alpha)=\int_{
\reel^{\prod_{D\in {\cal D}}\prod_{\gamma\in D}(|I_{\gamma}|-1)}}\pi_G(\theta|s,\alpha)
\;d\theta
\end{equation}
and where, as usual, ${\cal D}$ is equal to ${\cal E}$  when the model is saturated, to the set of complete subsets of $G$ when the model is  graphical Markov with respect to $G$ and  to the generating set for the model when the model is  hierarchical.

\noindent In order to be able to use this prior in practice, we need to answer a number of questions. The first basic question is to know for which values of the hyper parameters $(s, \alpha)$ where 
$$s=(s(i_D),\;D\in {\cal D}, i_D\in {\cal I}_D^*)\;\in\; \reel^{\prod_{D\in {\cal D}}\prod_{\gamma\in D}(|I_{\gamma}|-1)})\;\;\;\mbox{and}\;\;\;\alpha\in \reel $$ the distribution is proper, i.e. when does $I_G(s,\alpha)<+\infty$ hold. We will now give a necessary and sufficient condition for (\ref{pitheta}) to be proper as well as two practical methods to construct hyper parameters $(s,\alpha)$ such that it is proper. The next set of questions is concerned with the properties of this prior distribution in practice such as ease of prior specification, hyper Markov property. These questions will be addressed in \S5.

\subsection{A necessary and sufficient condition for the prior to be proper}
\begin{lemma}
\label{pigdetaillemma}
The prior distribution (\ref{pitheta}) is properif and only if  $(s,\alpha)$  belongs to 
\begin{equation}
\label{pigdetail}
\Pi_G=\{(s,\alpha)\;|\; \alpha>0, (\frac{s(i_D)}{\alpha}=\sum_{j\in {\cal I}, j_D=i_D}p(j),\; D\in {\cal D},\; i_D\in {\cal I}^*_D))\;\mbox{with} \;p(j)\;\mbox{as in (\ref{pe})} \}\;.
\end{equation}
\end{lemma}

\begin{pff}
Since the parameter space of (\ref{nefig}) is $\Theta_G=\reel^{\prod_{D\in {\cal D}}\prod_{\gamma\in D}(|{\cal I}_{\gamma}|-1)} $, by Theorem 1 of Diaconis and Ylvisaker (1974), a necessary and sufficient condition for (\ref{i}) to be  finite is that
$\alpha>0$ and $\frac{n}{\alpha}s=\frac{n}{\alpha}(s(i_D),\; D\in {\cal D},\; i_D\in {\cal I}_D)$ is in the interior of the convex hull of the support of $\mu_G$. Since the Laplace transform  
$$L_{\mu_G}(\theta)=\Big(1+\sum_{i \in {\cal I}, i\not = i^*}\exp \sum_{D\in {\cal D}^*_i}\theta_D(i_D)\Big)^n$$ 
is defined for $\Theta_G$
which is an open set, the  interior of the convex hull of the support of $\mu_G$ is equal to the mean space $M_G$ of ${\cal F}_{\mu_G}$. We therefore want to identify $M_G$. Let $k_{\mu_G}(\theta)=\log L_{\mu_G}(\theta)$.
 Since ${\cal F}_{\mu_G}$ is a natural exponential family with parameter 
 $\theta\in \Theta_G$, we have
\begin{eqnarray}
\label{diffm}
\hspace{-5mm}M_G=\{m=(m(i_D), \; D\in {\cal D},\; i_D\in {\cal I}^*_D)\;|\; m(i_D)&=&E(n(i_D))=n\frac{d k_{\mu_G}(\theta)}{d \theta(i_D)}=n\sum_{j\in {\cal I}, j_D=i_D}p(j)\;\}\;.
\end{eqnarray}
where $p(j)$ is as in (\ref{pe}). It follows immediately that $(s,\alpha)\in \Pi_G$ is a necessary and sufficient condition for  $\pi_G(\theta|s,\alpha)$ to be proper.
\end{pff}

\noindent From the lemma above, it is clear that in order to belong to $\Pi_G,\; (s,\alpha)$ must satisfy
$$\alpha\geq \mbox{max}_{D\in {\cal D}}s_D\;\;\;\mbox{and}\;\;s_D> s_E \;\;\mbox{for}\;\;D\subset E,\; D,E\in {\cal D}\;.$$
However, this condition is not sufficient since $(s,\alpha)$ must also be such that the $p(j)$ in (\ref{pigdetail}) satisfy the conditions $\theta_E(j_E)=0,\; E\not \in {\cal D}$.

\subsection{Two methods to construct $(s,\alpha)\in \Pi_G$}
From (\ref{pigdetail}), we immediately obtain the following method to construct hyper parameters $(s,\alpha)$ which are in $\Pi_G$:
\begin{enumerate}
\item[1.] Choose an arbitrary $\theta=(\theta(i_D),\; D\in {\cal D}|\; i_D\in {\cal I}_D^*)$
\item[2.] Compute $p(i)$ according to (\ref{pigeneralg}). 
\item[3.] Compute $\frac{s(i_D)}{\alpha}=\sum_{j\in {\cal I}, j_D=i_D}p(j)$ for $D\in {\cal D},\; i_D\in {\cal I}_D^*$.
\item[4.] Take $\alpha=1$.
\end{enumerate}
\vspace{2mm}

\noindent Another practical way to construct $(s,\alpha)\in \Pi_G$ is to start with a " prior contingency table" with all cell counts $n(i)$ positive. With $n$ denoting the total count in the given contingency table, the maximum likelihood estimate $\hat{p}$ of $p$ satisfying the equations
$$n(i_D)=n\sum_{j\in {\cal I}, j_D=i_D}\hat{p}(j), \; D\in {\cal D},\; i_D\in {\cal I}^*_D$$
and the constraints of the model, exists  and therefore we can take
$$\alpha=n,\;\; s(i_D)=n(i_D), \; D\in {\cal D},\; i_D\in {\cal I}^*_D\;$$
thus obtaining hyperparameters in $\Pi_G$. 

We note here that these hyperparameters are consistent across models since the "marginal counts" do not change when we take different models. Marginal counts do not change either when we take marginal or conditional models.

\section{The induced prior on the cell probabilities}
In this section, we will give the expression of the induced conjugate prior in terms of $p$, the cell probability parameter, first for graphical models Markov with respect to a decomposable $G$  thus making the link between the hyper Dirichlet and our conjugate prior, then for models Markov with respect to an arbitrary graph $G$  and finally for  general hierarchical models. 
\subsection{The conjugate prior when $G$ is decomposable}
For $G$ decomposable with set of cliques ${\cal C}=\{C_l,l=1,\ldots,k\}$ and set of minimal separators ${\cal S}=\{S_l,i=2,\ldots,k\}$, Dawid and Lauritzen (1993) defined the conjugate prior in terms of cell probabilities and called it the hyper Dirichlet distribution. Its density is expressed  in terms of 
\begin{equation}
\label{littlep}
p^{C_l}(i_D,i^*_{D^c}),\; D\subseteq C_l, l=1,\ldots,k,\;\;\;\;p^{S_l}(i_D,i^*_{D^c}),\; D\subseteq S_l,\; l=2,\ldots,k,\;\;\; D\in {\cal D}\;, i_D\in {\cal I}_D^*\;,
\end{equation}
the cell probabilities for the cliques and separators marginal tables, respectively. Note that in this subsection, for $D\subseteq C_l$ or $D\subseteq S_l, D^c$ denotes the complement of $D$ in $C_l$ or $S_l$ respectively. The density of the hyper Dirichlet is equal to
\begin{eqnarray}
&&\frac{\prod_{l=1}^k
\mbox{Dir}_{C_l}(p_{\emptyset}^{C_l}, p^{C_l}(i_D,i^*_{D^c});\alpha_{\emptyset}^{C_l},\alpha^{C_l}(i_D,i^*_{D^c}), D\in {\cal D}^{C_l}, i_D\in {\cal I}^*_D)
}{\prod_{l=2}^k \mbox{Dir}_{S_l}(p_{\emptyset}^{S_l},p^{S_l}(i_D,i^*_{D^c});\alpha^{S_l}(i_D,i^*_{D^c}), D\in {\cal D}^{S_l},i_D\in {\cal I}^*_D)
}\label{hyperdir}
\end{eqnarray}
with
\begin{eqnarray}
&&\mbox{Dir}_{C_l}(p_{\emptyset}^{C_l}, p^{C_l}(i_D,i^*_{D^c});\alpha^{C_l}(i_D,i^*_{D^c}), D\in {\cal D}^{C_l}, i_D\in {\cal I}^*_D)\nonumber\\
&=&\frac{\Gamma(\alpha_{\emptyset}^{C_l}+\sum_{D\in {\cal D}^{C_l}}\sum_{i_D\in {\cal I}^*_D}\alpha^{C_l}(i_D,i^*_{D^c}))}{\Gamma(\alpha_{\emptyset}^{C_l})\prod_{D\in {\cal D}^{C_l},i_D\in {\cal I}^*_D}\Gamma(\alpha^{C_l}(i_D,i^*_{D^c}))}(p_{\emptyset}^{C_l})^{\alpha_{\emptyset}^{C_l}-1}\prod_{D\in {\cal D}^{C_l}, i_D\in {\cal I}^*_D}(p^{C_l}(i_D, i^*_{D^c}))^{\alpha^{C_l}(i_D,i^*_{D^c})-1}\nonumber
\end{eqnarray}
with a similar expression for $\mbox{Dir}_{S_l}$ and where the hyper parameters
\begin{equation}
\label{hyperparametersdec}
(\alpha^{C_l}(i_D,i^*_{D^c}), D\in {\cal D}^{C_l}, i_D\in {\cal I}^*_D)\;\;\mbox{and} \;\;(\alpha^{S_l}(i_D,i^*_{D^c}), D\in {\cal D}^{S_l},i_D\in {\cal I}^*_D)
\end{equation}
 are hyperconsistent.

Since $\pi_G(\theta|s,\alpha)$ in (\ref{pitheta}) is the conjugate prior to the multinomial Markov with respect to $G$, it must coincide with the hyper Dirichlet when $G$ is decomposable. The aim of this subsection is to give the correspondence between the parameter $(s,\alpha)$ and the parameters of the hyper Dirichlet  explicitly.

The probabilities in (\ref{littlep}) are not all free variables since, by the Markov properties of the multinomial distribution,
$$p(i)=\frac{\prod_{l=1}^kp^{C_l}(i_{C_l})}{\prod_{l=2}^kp^{S_l}(i_{S_l})}\;$$ 
 and therefore some are functions of the others. Let ${\cal D}^{C_i}$ and ${\cal D}^{S_i}$ denote the set of nonempty subsets of $C_i$ and $S_i$ respectively. We can choose the free marginal probabilities to be 
\begin{equation}
\label{chosenpdiscrete}
p^G=(p^{C_l}(i_D,i^*_{D^c}), D\in {\cal D}^{C_l}\setminus \cup_{j=2}^k {\cal D}^{S_j},l=1,\ldots,k, p^{S_l}(i_D,i^*_{D^c}), D\in {\cal D}^{S_l}, l=2,\ldots,k,i_D\in {\cal I}_D^*)\;.
\end{equation}

\noindent The Jacobian of the change of variable $\theta\mapsto p^G$ is given in the following lemma.

\begin{lemma}
\label{lemmajacobian}
The Jacobian of the change of variables from $\theta=(\theta(i_D), D\in {\cal D}, i_D\in {\cal I}_D^*)$  as given in (\ref{pitotheta}) to  $p^G$ as given in (\ref{chosenpdiscrete}) is
\begin{equation}
\label{hyperJdiscrete}
\Big|\frac{d \theta}{d p^G}\Big|^{-1}=\frac{\prod_{l=1}^kp^{C_l}_{\emptyset}\prod_{ D\in {\cal D}^{C_l}}\prod_{i_D\in {\cal I}_D^*}p^{C_l}(i_D)}
{\prod_{l=2}^kp^{S_l}_{\emptyset}\prod_{D\in {\cal D}^{S_l}}\prod_{i_D\in {\cal I}_D^*}p^{S_l}(i_D)}
\end{equation}
\end{lemma}
The proof of this lemma is given in the Appendix.
The correspondence between $(s,\alpha)$ and (\ref{hyperparametersdec}) is given by the following proposition.

\begin{prop}
\label{alphadec}
When the graph $G$ is decomposable with set of cliques $(C_i,i=1,\ldots,k)$ and sets of minimal separators $(S_i, i=2,\ldots,k)$, the conjugate prior induced from (\ref{pitheta}) is identical to the hyper Dirichlet (\ref{hyperdir}) with hyper parameters (\ref{hyperparametersdec})
where
\begin{eqnarray}
\hspace{-5mm}\alpha^{C_l}(i_D,i^*_{D^c})=\sum_{C_l\supseteq F\supseteq D}\sum_{j_F\in {\cal I}_F^*\;|\;(j_F)_D=i_D}(-1)^{|F\setminus D|}s(j_F)
&&\alpha_{\emptyset}^{C_i}=\alpha+\sum_{D\subseteq C_i}(-1)^{|D|}\sum_{i\in {\cal I}_D^*}s(i_D)\label{adci}\\
\hspace{-5mm}\alpha^{S_l}(i_D,i^*_{D^c})=\sum_{S_l\supseteq F\supseteq D}\sum_{j_F\in {\cal I}_F^*\;|\;(j_F)_D=i_D}(-1)^{|F\setminus D|}s(j_F)&&\alpha_{\emptyset}^{S_l}=\alpha+\sum_{D\subseteq S_l}(-1)^{|D|}\sum_{i\in {\cal I}_D^*}s(i_D)\label{adsi}\;.
\end{eqnarray}
Moreover
\begin{equation}
I_G(s,\alpha)=\frac{
\prod_{l=1}^k \Gamma(\alpha_{\emptyset}^{C_l})
\prod_{ D\in {\cal D}^{C_l}}\prod_{i_D\in {\cal I}^*_D}   \Gamma(\alpha^{C_l}(i_D,i^*_{D^c}))
}
{\prod_{l=2}^k\Gamma(\alpha_{\emptyset}^{S_l})
\prod_{ D\in {\cal D}^{S_l}}\prod_{i_D\in {\cal I}^*_D}   \Gamma(\alpha^{S_l}(i_D,i^*_{D^c}))
}\label{iratio}
\end{equation}
\end{prop}
\begin{pff}
\noindent Since the distribution of $Y$ in (\ref{nefmarkov}) is Markov with respect to $G$, we have that
\begin{equation}
\label{pdmarkov}
p(i)=\frac{\prod_{l=1}^k p^{C_l}(i_{ C_l})}{\prod_{l=2}^k p^{S_l}(i_{S_l})}\;.
\end{equation}
Then
\begin{eqnarray}
\theta_E(i_E)&=&\sum_{F\subseteq E}(-1)^{|E\setminus F|}\log p(i_F,i^*_{F^c})\nonumber\\
&=&\sum_{F\subseteq E}(-1)^{|E\setminus F|}
\Big(\sum_{l=1}^k \log p^{C_l}(i_{F\cap C_l},i^*_{F^c\cap C_l})-\sum_{l=2}^k \log p^{S_l}(i_{F\cap S_l},i^*_{F^c\cap S_l})\Big)\nonumber\\
&=&\sum_{l=1}^k\Big(\sum_{F\subseteq E}(-1)^{|E\setminus F|}\log p^{C_l}(i_{F\cap C_l},i^*_{F^c\cap C_l})\Big)-
\sum_{l=2}^k\Big(\sum_{F\subseteq E}(-1)^{|E\setminus F|}\log p^{S_l}(i_{F\cap S_l},i^*_{F^c\cap S_l})\Big)\nonumber\\
&=&\sum_{l=1}^k\theta^{C_l}_E(i_{E\cap C_l})-\sum_{l=2}^k\theta^{S_l}_E(i_{E\cap S_l})\nonumber
\end{eqnarray}
If $E\subseteq C_l$, $E\cap C_l=E$ and $\theta^{C_l}_E(i_{E\cap C_l})=\theta^{C_l}_E$. If $E\not \subseteq C_l$, then by Lemma \ref{glim}, $\theta^{C_l}_E(i_{E\cap C_l})=0$ and similarly for $\theta^{S_l}_E(i_{E\cap S_l})$. We therefore have

\begin{equation}
\label{decomptheta}
\theta_D(i_D)=\sum_{l=1}^k \theta_D^{C_i}(i_D)-\sum_{l=2}^k \theta_D^{S_l}(i_D),\;D\in {\cal D}
\end{equation}
where 
\begin{eqnarray}
\theta^{C_i}_D(i_D)&=& \sum_{F\subseteq D, F\in {\cal D}_0^{C_i}}  (-1)^{|D\setminus F|}  \log p^{C_i}(i_F, i^*_{F^c}),\;\;\;\mbox{for}\;\;D\subseteq C_i\nonumber\\
\theta_{\emptyset}^{C_i}&=&\log p_{\emptyset}^{C_i}\nonumber\\
\theta^{C_i}_D(i_D)&=&0\;\;\;\mbox{for}\;\;D\not \subseteq C_i\nonumber
\end{eqnarray}
and similar expressions for $\theta^{S_i}(i_D)$ (see also Consonni and Leucari, 2005 for the derivation of these formulas in the case of bivariate data). From (\ref{pdmarkov}), we also have
\begin{eqnarray}
\log \Big(1+\sum_{j \in {\cal I}, j\not = i^*}\exp\sum_{D\in {\cal D}^*_j}\theta_D(j_D)\Big)&=&-\log p_{\emptyset}=-\Big(\sum_{l=1}^k \log p_{\emptyset}^{C_l}-\sum_{l=2}^k \log p_{\emptyset}^{S_l}\Big)
\end{eqnarray}
Therefore  (\ref{pitheta}) can be written as

\begin{eqnarray}
\pi_{G}(\theta(p)|s,\alpha)&\propto&\frac{\prod_{l=1}^k\exp\Big\{\sum_{D\in {\cal D}^{C_l}}\sum_{i_D\in {\cal I}_D^*}\Big(\sum_{E\subseteq D}(-1)^{|D\setminus E|}\log p^{C_l}(i_E, i^*_{E^c})\Big)s(i_D)+\alpha\log p^{C_l}_{\emptyset}\Big\}}
{\prod_{l=2}^k\exp\{\sum_{D\in {\cal D}^{S_l}}\sum_{i_D\in {\cal I}_D^*}\Big(\sum_{E\subseteq D}(-1)^{|D\setminus E|}\log p^{S_l}(i_E, i^*_{E^c})\Big)s(i_D)+\alpha\log p^{S_l}_{\emptyset}\}
}\nonumber\\
&=&\frac{\prod_{l=1}^k\exp\Big\{\sum_{E\in {\cal D}^{C_l}}\sum_{i_E\in {\cal I}_E^*}\alpha^{C_l}(i_E,i^*_{E^c})\log p^{C_l}(i_E, i^*_{E^c})
+\alpha^{C_l}_{\emptyset}\log p^{C_l}_{\emptyset}\Big\}}
{\prod_{l=2}^k\exp\Big\{\sum_{E\in {\cal D}^{S_l}}\sum_{i_E\in {\cal I}_E^*}\alpha^{S_l}(i_E,i^*_{E^c})\log p^{S_l}(i_E, i^*_{E^c})
+\alpha^{S_l}_{\emptyset}\log p^{S_l}_{\emptyset}\Big\}}\nonumber\\
\label{pithetaofp}
&=&\frac{\prod_{l=1}^k(p^{C_l}_{\emptyset})^{\alpha^{C_l}_{\emptyset}}\prod_{E\in {\cal D}^{C_l}}\prod_{i_E\in {\cal I}_E^*}\;\;(p^{C_l}(i_E, i^*_{E^c}))^{\alpha^{C_l}}(i_E,i^*_{E^c})}
{\prod_{l=2}^k(p^{S_l}_{\emptyset})^{\alpha^{S_l}_{\emptyset}}\prod_{E\in {\cal D}^{S_l}}\prod_{i_E\in {\cal I}_E^*}\;\;(p^{S_l}(i_E, i^*_{E^c}))^{\alpha^{S_l}}(i_E,i^*_{E^c})}
\end{eqnarray}
where $\alpha^{C_l}(i_E,i^*_{E^c}),\;\alpha^{S_l}(i_E,i^*_{E^c}),\;\alpha^{C_l}_{\emptyset}$ and $\alpha^{S_l}_{\emptyset}$ are as defined in (\ref{adci}) and (\ref{adsi}).

\noindent The induced prior on $p$ is obtained by multiplying (\ref{pithetaofp}) by the Jacobian
 (\ref{hyperJ}) and it follows immediately that it is the hyper Dirichlet with
hyper parameters as given in (\ref{adci}) and (\ref{adsi}).

The expression of (\ref{iratio}) is obtained by noticing that for any $C_i$ or $S_i$,
$$\alpha^{C_l}_{\emptyset}+\sum_{E\in {\cal D}^{C_l}}\sum_{i_E\in {\cal I}^*_E}
\alpha^{C_i}(i_E,i^*_{E^c})=\alpha=\alpha^{S_l}_{\emptyset}+\sum_{E\in {\cal E}^{S_i}}\sum_{i_E\in {\cal I}^*_E}\alpha^{S_l}(i_E,i^*_{E^c})\;.$$
This completes the proof.
\end{pff}
\subsection{The conjugate prior when $G$ is arbitrary}

To obtain the conjugate prior in terms of $p$, we need to compute the Jacobian
$\frac{d \theta}{d p}$
of the transformation from $\theta$ to $p$ as defined respectively in (\ref{freethetadiscrete}) and (\ref{freepdiscrete}). Before doing so, we need to define the following quantities. For $C\in {\cal D}, \;H\in {\cal E}$, let
\begin{equation}
\label{fich}
F(i_C,j_H)=\left \{
\begin{array}{ll}(-1)^{|C|-1}& \mbox{if}\;\; (j_H)_C=i_C \\0 & \mbox{otherwise}\end{array} \;\;\;\;\mbox{for}\;\;i_C\in {\cal I}_C^*, \;j_H\in {\cal I}^*_H\;.\right .
\end{equation}
These $F(i_C,j_H)$ can be gathered in a $\prod_{D\in {\cal D}}|{\cal I}^*_D|\;\times \;\prod_{H\in {\cal E}}|{\cal I}^*_H|$ matrix $F$ where the rows are indexed by $i_D\in {\cal I}^*_D, D\in {\cal D}$ and the columns by $j_H\in {\cal I}^*_H, H\in {\cal E}\cup \{\emptyset\}$.
\vspace{2mm}

\noindent For example, in the case of binary data for ${\cal D}$ and ${\cal E}$ as given in \S 2.4
 the matrix $F$ is
\begin{equation}
\label{fexample}
F=\pmatrix{0&1&0&0&0& 1&0 &0 &1 &1&0&1 &0 &1 &1 &1 \cr
           0&0&1&0&0& 1&1 &0 &0 &0&1&1 &1 &0 &1 &1 \cr
           0&0&0&1&0& 0&1 &1 &0 &1&0&1 &1 &1 &0 &1 \cr
           0&0&0&0&1& 0&0 &1 &1 &0&1&0 &1 &1 &1 &1 \cr
           0&0&0&0&0&-1&0 &0 &0 &0&0&-1&0 &0 &-1&-1\cr
           0&0&0&0&0& 0&-1&0 &0 &0&0&-1&-1&0 &0 &-1\cr
           0&0&0&0&0& 0&0 &-1&0 &0&0&0 &-1&-1&0 &-1\cr
           0&0&0&0&0& 0&0 &0 &-1&0&0&0 &0 &-1&-1&-1} 
\end{equation}
\vspace{2mm}

We also need the following two lemmas. Their proof is given in the Appendix.
\begin{lemma}
\label{F}
Let $G$ be a nondecomposable prime graph. For the matrix $F$ as described  in (\ref{fich}), the sum of the entries in each column $j_H,\;j\in {\cal I}^*_H,\; H\in {\cal E}$ is such that
\begin{equation}
\label{sum}
\sum_{i_C\in {\cal I}^*_C, C\in {\cal D}}F(i_C,j_H)=1
\end{equation}
 if and only if $H$, as an induced subgraph of $G$, is decomposable and connected.
\end{lemma}

We are now in a position to give the expression of the Jacobian. Let 
$${\cal U}=\{F\in {\cal E}\;| F\;\;\mbox{is either nondecomposable or nonconnected}\;\}$$
and we also write ${\cal U}_{0}={\cal U}\cup \{\emptyset\}.$
\begin{lemma}
\label{jacob}
Let 
\begin{equation}
\label{ah}
a(j_H)=\Big(\sum_{i_C\in {\cal I}^*_C, C\in {\cal D}}F(i_C,j_H)-1\Big),\;\;j\in {\cal I}^*_H,\;H\in {\cal E}\cup \emptyset\;.
\end{equation}
The Jacobian of the transformation
\begin{equation}\label{changeptheta}
p=(p(i_D,i^*_{D^c}), \; i_D\in {\cal I}^*_D,\;D\in {\cal D})\mapsto \;\theta=(\theta_D(i_D),\; i_D\in {\cal I}^*_D, \; D\in {\cal D}),
\end{equation}
where $p$ is as given in (\ref{freepdiscrete}) and $\theta$ as in (\ref{pitotheta}) is
\begin{eqnarray}
J=\Big|\frac{d p}{d \theta}\Big|&=&\Big(\prod_{D\in {\cal D}}\prod_{i_D\in {\cal I}^*_D}p(i_D,i^*_{D^c})\Big)\Big(-\sum_{H\in {\cal U}_0}\sum_{j_H\in {\cal I}_H^*}a(j_H,j^*_{H^c})\frac{p(j_H,j^*_{H^c})}{p_{\emptyset}}\Big)\nonumber\\
&=&\Big(p_{\emptyset}\prod_{D\in {\cal D}}\prod_{i_D\in {\cal I}^*_D}p(i_D,i^*_{D^c})\Big)\Big(1-\frac{1}{p_{\emptyset}}\sum_{H\in {\cal U}}
\sum_{j_H\in {\cal I}_H^*}a(j_H,j^*_{H^c})p(j_H,j^*_{H^c})\Big)
\label{jacobianptheta}
\end{eqnarray}
\end{lemma}

\noindent The proof of this lemma is given in the Appendix.
\vspace{2mm}

\noindent We can now give the conjugate prior (\ref{pitheta}) in terms of $p$ as given in (\ref{freepdiscrete}). 
Let us note first that by (\ref{pitotheta}) and (\ref{nefig}), the marginal cell counts $y=(n(i_D), i_D\in {\cal I}_D^*,\;D\in {\cal D})$ for the multinomial distribution Markov with respect to $G$ has density 

\begin{eqnarray}
f((y\;|\;p)&\propto&\prod_{D\in {\cal D}}\prod_{i_D\in {\cal I}_D^*}\Big(\prod_{F\subseteq D}p(i_F, i^*_{F^c})^{(-1)^{|D\setminus F|}}\Big)^{y(i_D)}p_{\emptyset}^{n}\nonumber
\\
&\propto&\prod_{D\in {\cal D}}\prod_{i_D\in {\cal I}_D^*}
p(i_D,i^*_{D^c})^{u(i_D)}
p_{\emptyset}^{u(\emptyset)}\label{multid}
\end{eqnarray}
where $u(i_D)=\sum_{F\supseteq D}(-1)^{F\setminus D}\sum_{j_F\;|\;(j_F)_D=i_D}
y(j_F)$ and $u(\emptyset)=n+\sum_{i_D\in {\cal I}^*_D, D\in {\cal D}}(-1)^{|D|}y(i_D)$.
\begin{theorem}
\label{inducedprior}
For $(s,\alpha)\in \Pi_G$ and $a_H$ as given in (\ref{ah}), the conjugate prior distribution induced from (\ref{pitheta}) by (\ref{changeptheta}), that is, the conjugate prior for the  parameter $p$ of the multinomial family of distributions  (\ref{multid}) is
\begin{eqnarray}
\label{piptheta}
\pi^p_G(p\;|\; (s,\alpha))&=&\frac{K^{-1}}{I_G(s,\alpha)}\prod_{D\in {\cal D}}\prod_{i_D\in {\cal I}^*_D}p(i_D,i^*_{D^c})^{\alpha(i_D,i^*_{D^c})
-1}p_{\emptyset}^{\alpha_{\emptyset}-1}\label{pip}\;.
\end{eqnarray}
where  
\begin{eqnarray}
K&=&\Big(1-\frac{1}{p_{\emptyset}}\sum_{H\in {\cal U}}
\sum_{j_H\in {\cal I}_H^*}a(j_H,j^*_{H^c})p(j_H,j^*_{H^c})\Big)\nonumber\\
\alpha(i_D,i^*_{D^c})&=&\sum_{F\supseteq D}\sum_{j_F\in {\cal I}_F^*\;|\;(j_F)_D=i_D}(-1)^{F\setminus D}s(i_F)\label{alphaididc}\\
\alpha_{\emptyset}&=&\alpha+\sum_{i_D\in {\cal I}^*_D, D\in {\cal D}}(-1)^{|D|}s(i_D)\;.\label{alphanot}
\end{eqnarray}
\end{theorem}
This result follows immediately from the expression of the conjugate prior (\ref{pitheta}) in terms of $\theta$, (\ref{pitotheta}) and (\ref{jacobianptheta}).

\vspace{2mm}

\noindent {\bf Example}

\noindent When the graph is the four cycle with binary data as considered before
$${\cal U}_0=\{ac, bd, abcd, \emptyset   \}.$$
From (\ref{fexample}) and the constraints $\theta_E(i_E)=0$ for $E\not\in {\cal D}$, we have
$$a_{ac}=a_{bd}=1,\;\; a_{abcd}=-1,\;\;\frac{p_{ac}}{p_{\emptyset}}=\frac{p_ap_c}{p_{\emptyset}^2}, \;\;\frac{p_{bd}}{p_{\emptyset}}=\frac{p_bp_d}{p_{\emptyset}^2},\;\;\frac{p_{abcd}}{p_{\emptyset}}=\frac{p_{ab}p_{bc}p_{cd}p_{da}}{p_ap_bp_cp_d}$$
and
\begin{eqnarray}
\lefteqn{\pi(p_D, D\in {\cal D}\;|\; (s,\alpha))=}\nonumber\\
&&I_G(s,\alpha)^{-1}p_a^{s_a-s_{da}-s_{ab}-1}p_b^{s_b-s_{ab}-s_{bc}-1}p_c^{s_c-s_{bc}-s_{cd}-1}p_d^{s_d-s_{cd}-s_{da}-1}p_{ab}^{s_{ab}-1}p_{bc}^{s_{bc}-1}p_{cd}^{s_{cd}-1}p_{da}^{s_{da}-1}p_{\emptyset}^{\alpha-1}\nonumber\\
&&\hspace{1cm}(1-\frac{p_ap_c}{p_{\emptyset}^2}-\frac{p_bp_d}{p_{\emptyset}^2}+
\frac{p_{ab}p_{bc}p_{cd}p_{da}}{p_ap_bp_cp_d})^{-1}\nonumber
\end{eqnarray}
\subsection{The conjugate prior for a general hierarchical model}
When the model is not specified to be graphical but is a hierarchical log-linear model, we can also obtain the induced prior in terms of $p$ and the statement is similar to Theorem \ref{inducedprior} above except that the term coming from the Jacobian $|\frac{d \theta}{d p}|$ is more general and we have

\begin{theorem}
\label{inducedpriorhierarchical}
For $(s,\alpha)\in \Pi_G$  the conjugate prior distribution induced from (\ref{pitheta}) by (\ref{changeptheta}), that is, the conjugate prior for the  parameter $p$ of the multinomial family of distributions  (\ref{multid}) for the hierarchical log-linear model is as in (\ref{piptheta}) with
 
\begin{eqnarray}
K&=&1-\sum_{H\in {\cal E}}\sum_{j_H\in {\cal I}_H^*}p(j_H,j^*_{H^c})(\sum_{\{D\subseteq H, D\in {\cal D}\}}\sum_{\{C\subseteq D, C\in {\cal D}\}}(-1)^{|D\setminus C|})\nonumber
\end{eqnarray}
and $\alpha(i_D,i^*_{D^c})$ and $\alpha_{\emptyset}$ as in (\ref{alphaididc}) and (\ref{alphanot})
\end{theorem}
The proof follows immediately from the expression of the conjugate prior (\ref{pitheta}) in terms of $\theta$, (\ref{pitotheta}), (\ref{jacobianptheta}) and Remark (\ref{notgraphical}).

\section{Properties of the conjugate prior}

\subsection{Hyper-parameter specification}
Let us now turn to the practical problem of choosing hyperparameters which will reflect either some prior belief or lack of prior belief.

Suppose first that we do not have any prior information and want to put a flat prior on the log-linear parameters. From the expression (\ref{fg}) of the distribution of the marginal counts $y=(n(i_D), D\in {\cal D}, i_D\in {\cal I}_D^*)$, it is clear that the hyperparameters $s(i_D)$ can be thought of as the "prior marginal counts" for the marginal cell $i_D$. Therefore, we can take for $s(i_D)$ the set of $i_D$-"marginal counts" ,  $D\in {\cal D}, i_D\in {\cal I}_D^*$ for a "prior" contingency table with all "cell counts"' equal to $\frac{1}{\cal I}$. We can also take $\alpha$ to be the "total count", that is 1 in this case. This would lead, of course to
$$s(i_D)=\sum_{j\in {\cal I}, j_D=i_D\in {\cal I}_D^*}\frac{1}{\cal I}=\frac{1}{\cal I}\prod_{\gamma \in D^c}|{\cal I}_{\gamma}|\;,$$
where $D^c$ is the complement of $D$ in $V$.
Since for the saturated model, the conjugate prior in terms of cell probabilities is the Dirichlet, it is clear that this choice of hyperparameters also yields a flat prior for the cell probabilities of the saturated model with all hyperparameters being equal to $\frac{1}{\cal I}$. This prior is in fact the vague prior advocated by Perks (1947) (see also Dellaportas and Forster, 1999).

If we have prior information, we can first exclude all the interactions that are thought to be absent. Indeed, by Lemma \ref{hamcli}, if two variables are believed to be independent given the others, then all $\theta_E(i_E)=0$ for $E\in {\cal E}$ containing these two variables. We may have additional information such as  the knowledge of positive or negative interaction between one or more variables. This knowledge can be expressed by computing the expected value and variance for $e^{\theta_D}$ for appropriate $D\in {\cal D}$.
 To illustrate what we mean, let us consider the data given by Hook, Albright and Cross (1980)  and studied by King and Brooks (2001). In this data set, there are three variables 
$$a\equiv BC,\;\;b\equiv DC,\;\;c\equiv MR$$ each taking the values 1 or 0 representing the presence or absence of, respectively, birth certificates, death certificates and medical records for each individual. The individuals under study are children with spina bifida. The data  consists of an incomplete contingency table for each one of  six  years. From Hook, Albright and Cross (1980), it can reasonably be assumed that the model is the decomposable graphical model with cliques $a$ and $bc$. Since the data is binary, from (\ref{nefmarkov}),  the conjugate prior will then be of the form
\begin{eqnarray}
\label{pigexample}
\pi_G(\theta|s,\alpha)&=&I_G(s,\alpha)^{-1}\exp\;\Big(\theta_as_a+\theta_bs_b+\theta_cs_c+\theta_{bc}s_{bc}\\
&&\hspace{1cm}-\alpha\log (1+e^{\theta_a}+e^{\theta_b}+e^{\theta_c}+e^{\theta_a+\theta_b}+e^{\theta_a+\theta_c}+e^{\theta_b+\theta_c+\theta_{bc}}+e^{\theta_a+\theta_b+\theta_c+\theta_{bc}}\Big)\nonumber
\end{eqnarray}
There is also some prior knowledge about the interaction between $b$ and $c$, that is for 
$$e^{\theta_{bc}}=\frac{p_{bc}p_{\emptyset}}{p_bp_c}.$$
With high probability, $e^{\theta_{bc}}$ is expected to be in the interval $(-.9,-.1)$.
From (\ref{pigexample}) and the formulas given  in Proposition \ref{alphadec}, if we let $s'=(s_a,s_b,s_c,s_{bc}+1)$ we have
\begin{eqnarray}
E(e^{\theta_{bc}})&=&\frac{I_G(s',\alpha)}{I_G(s,\alpha)}\nonumber\\
&=&\frac{\Gamma(s_{bc}+1)\Gamma(s_b-s_{bc}-1)\Gamma(s_c-s_{bc}-1)\Gamma(\alpha-s_b-s_c+s_{bc}+1)}
{\Gamma(s_{bc})\Gamma(s_b-s_{bc})\Gamma(s_c-s_{bc})\Gamma(\alpha-s_b-s_c+s_{bc})}\nonumber\\
&=&\frac{s_{bc}(\alpha-s_b-s_c+s_{bc})}{(s_b-s_{bc}-1)(s_c-s_{bc}-1)}\label{ethetabc}
\;.
\end{eqnarray}
We therefore have the constraint
$$-.9\leq \frac{s_{bc}(\alpha-s_b-s_c+s_{bc})}{(s_b-s_{bc}-1)(s_c-s_{bc}-1)}\leq -.1.$$

In the absence of any prior knowledge on the other log-linear parameters, we can assume that their expectation is around
0 which would imply that
\begin{eqnarray}
E(e^{\theta_{a}})&=&\frac{\Gamma(\alpha-s_a-1)\Gamma(s_a+1)}{\Gamma (\alpha-s_a)\Gamma(s_a)}=              \frac{s_a}{\alpha-s_a-1}\nonumber\\
E(e^{\theta_{b}})&=&\frac{\Gamma(\alpha-s_b-1-s_c+s_{bc})\Gamma(s_b+1-s_{bc})}{\Gamma (\alpha-s_b-s_c+s_{bc})\Gamma(s_b-s_{bc})}=\frac{s_b-s_{bc}}{\alpha-s_b-s_c+s_{bc}-1}\nonumber\\
E(e^{\theta_{c}})&=&\frac{\Gamma(\alpha-s_b-1-s_c+s_{bc})\Gamma(s_c+1-s_{bc})}{\Gamma (\alpha-s_b-s_c+s_{bc})\Gamma(s_c-s_{bc})}=\frac{s_c-s_{bc}}{\alpha-s_b-s_c+s_{bc}-1}\nonumber
\end{eqnarray}
are all around 1.
If we took all three ratios to be 1, we would obtain the relationships
\begin{eqnarray}
&&2s_a=\alpha-1,\;s_b-s_{bc}=s_c-s_{bc},\; 2(s_b-s_{bc})=\alpha-s_c-1,\; 2(s_c-s_{bc})=\alpha-s_b-1,\nonumber\\
&& -.9\leq \frac{s_{bc}}{s_b-s_{bc}-1}\leq -.1\;.
\end{eqnarray}
and choose appropriate $(s,\alpha)$ satisfying these conditions. We might also want to compute the variance of these quantities which is, of course, also immediate with the results of Proposition \ref{alphadec}, and give an interval where we wish $E(e^{\theta_D}), D\in \{a,b,c\}$ to be.

In general, when the model considered is not necessarily a decomposable graphical model, the ratio of normalising constants of the type $\frac{I_G(s',\alpha)}{I_G(s,\alpha)}$ has to be computed numerically. This is feasible by any or the standard MCMC or approximation methods. However, it might be wiser and much simpler to choose a decomposable model covering the interaction believed to be true. For example, if, in the example above, the prior model was believed to be the hierarchical model with generating class $\{ab, bc,ca\}$, then a reasonable prior model would be the saturated model Markov with respect to the complete graph subject to the fact that the interaction between $a,b$ and $c$ is weak, that is $E(e^{\theta_{abc}})$ is close to 1.

It  remains to know whether the hyperparameters chosen for the conjugate prior on the log-linear parameter $\theta$ will yield hyper parameters in the conjugate prior induced by (\ref{pitheta}) for the cell probabilities which are consistent with the given prior beliefs. From Theorem \ref{inducedprior}, we know that the induced prior for the cell probabilities "`looks" like a Dirichlet on the free cell probabilities, that is $p(i_D, i^*_{D^c}), D\in {\cal D}, i\in {\cal I}^*_D$ with an additional factor for the Jacobian.
The powers of the $p(i_D, i^*_{D^c})$ correspond to "prior cell counts" $n(i_D)$ and therefore any choice $s(i_D)$ in (\ref{pitheta}) will have the same meaning in  (\ref{piptheta}). For example, corresponding to the condition that (\ref{ethetabc}) be in the interval $(-.9,-.1)$ corresponds the condition that
$$E(\frac{p_{bc}p_{\emptyset}}{p_bp_c})$$ be in that interval also. From (\ref{piptheta}),
the conjugate prior on $p=(p_a,p_b,p_c,p_{bc})$ is
\begin{eqnarray}
\pi^p_G(p\;|\; (s,\alpha))&=&\frac{\Big(1-\frac{p_ap_b+p_ap_c+p_ap_{bc}}{p_{\emptyset}^2}\Big)^{-1}}{I_G(s,\alpha)}
p_a^{s_a-1}p_b^{s_b-s_{bc}-1}p_c^{s_c-s_{bc}-1}p_{bc}^{s_{bc}-1}p_{\emptyset}^{\alpha-s_a-s_b-s_c+s_{bc}-1}\nonumber\;.
\end{eqnarray}
Therefore
$$E(\frac{p_{bc}p_{\emptyset}}{p_bp_c})=\frac{I_G(s',\alpha)}{I_G(s,\alpha)}$$
where $s'=(s_a,s_b-1,s_c-1,s_{bc}+1)$ and it follows immediately that
$$E(\frac{p_{bc}p_{\emptyset}}{p_bp_c})=\frac{s_{bc}(\alpha-s_b-s_c+s_{bc})}{(s_b-s_{bc}-1)(s_c-s_{bc}-1)}\;,$$
thus giving the same condition as in (\ref{ethetabc}).

\subsection{The strong hyper Markov property for local updates in graphical model}
Let us now assume that the multinomial distribution of the contingency cell counts  is Markov with respect to  an arbitrary undirected graph $G$. We know from Dawid and Lauritzen (1993) that the multinomial distribution is strong meta Markov and as the conjugate distribution of the parameter $\theta$ of the exponential family (\ref{nefmarkov}), the conjugate prior (\ref{pitheta}) is strong hyper Markov. 

Consider the decomposition of $G$ into its prime components and let $P_l, l=1,\ldots,k$ be a perfect enumeration of these components. Let $S_l, l=2,\ldots,k$ be the corresponding separators. We now give the expression of (\ref{pitheta}) as the Markov ratio of conjugate  priors on the prime components and the separators of $G$.
 
\begin{prop}
The conjugate prior (\ref{pitheta}) can be written as the Markov ratio
\begin{eqnarray}
\label{pithetap}
\pi_G(\theta|s,\alpha)&=&\frac{\prod_{l=1}^k \pi_{P_l}(\theta^{P_l}|s^{P_l},\alpha)
}{\prod_{l=2}^k\pi_{S_l}(\theta^{S_l}|s^{S_l},\alpha)}
\end{eqnarray}
where
\begin{eqnarray}
\label{pithetapl}
\lefteqn{ \pi_{P_l}(\theta^{P_l}|s^{P_l},\alpha)}\\
&=&I_G(s^{P_l},\alpha)^{-1}\exp \{\sum_{D\in {\cal D}_{P_l}}\sum_{i_D\in {\cal I}_D^*}\theta_D^{P_l}(i_D)s(i_D)\}-\alpha\log \Big(1+\sum_{i \in {\cal I}_{P_l}, i\not = i^*}\exp \sum_{D\in {\cal D}^*_i}\theta_D^{P_l}(i_D)\Big)\nonumber
\end{eqnarray}
 and where $s^{P_l}=(s(i_D),\;{\cal D}^{P_l}, i_D\in {\cal I}^*_D)$ and $s^{S_l}=(s(I_D),\;{\cal D}^{S_l}, i_D\in {\cal I}^*_D)$.

 \noindent The induced conjugate prior (\ref{piptheta}) on $p$  can be written as the corresponding Markov ratio of conjugate priors induced on $p^{P_l}$ and $p^{S_l}$ from $\pi_{P_l}(\theta^{P_l}|s^{P_l},\alpha)$ in (\ref{pithetapl}) and 
 $\pi_{S_l}(\theta^{S_l}|s^{S_l},\alpha)$.
  \end{prop}
\begin{pff}
It is not difficult to see that
\begin{equation}
\label{decompthetap}
\theta_D(i_D)=\sum_{l=1}^k \theta_D^{P_l}(i_D)-\sum_{l=2}^k \theta_D^{S_l}(i_D),\;D\in {\cal D}
\end{equation}
where 
\begin{eqnarray}
\theta^{C_l}_D(i_D)&=& \sum_{F\subseteq D, F\in {\cal D}_0^{C_l}}  (-1)^{|D\setminus F|}  \log p^{C_l}(i_F, i^*_{F^c}),\;\;\;\mbox{for}\;\;D\subseteq C_l\nonumber\\
\theta_{\emptyset}^{P_l}&=&\log p_{\emptyset}^{P_l}\nonumber\\
\theta^{P_l}_D(i_D)&=&0\;\;\;\mbox{for}\;\;D\not \subseteq P_l,\;l=1,\ldots,k,\;\;\;\;\;\theta^{S_l}_D(i_D)=0\;\;\;\mbox{for}\;\;D\not \subseteq S_l, l=2,\ldots,k, i_D\in {\cal I}^*_D\nonumber
\end{eqnarray}
We have proved this property for a decomposable graph $G$ in \S 4.1, (\ref{decomptheta}). The proof goes exactly along the same lines here. Therefore, if we let
$${\cal D}^{P_l}=\{D\in {\cal D}\;|\; D\subseteq P_l\},\;\;\;\;{\cal D}^{S_l}=\{D\in {\cal D}\;|\; D\subseteq S_l\}, \;$$
$$\theta^{P_l}=(\theta^{P_l}(i_D), D\in {\cal D}^{P_l}, i_D\in {\cal I}^*_D),\;\;\;\;
\theta^{S_l}=(\theta^{S_l}(i_D), D\in {\cal D}^{S_l}, i_D\in {\cal I}^*_D),$$
and
$$s^{P_l}=(s(i_D), D\in {\cal D}^{P_l}, i_D\in {\cal I}^*_D),\;\;\;\;
s^{S_l}=(s(i_D), D\in {\cal D}^{S_l}, i_D\in {\cal I}^*_D),$$
we see that
\begin{eqnarray}
\sum_{D\in {\cal D}, i_D\in {\cal I}^*_D}s(i_D)\theta_D(i_D)&=&\sum_{l=1}^k\sum_{D\in {\cal D}^{P_l}, i_D\in {\cal I}_D^*}\theta_D^{P_l}(i_D)s(i_D)-\sum_{l=2}^k\sum_{D\in {\cal D}^{S_l}, i_D\in {\cal I}_D^*}\theta_D^{S_l}(i_D)s(i_D)\nonumber\;.
\end{eqnarray}

Since by the Markov property, we also have 
$p_{\emptyset}=\frac{
\prod_{l=1}^k p_{\emptyset}^{P_l}
}{
\prod_{l=2}^k p_{\emptyset}^{S_l}
}$, that is

\begin{eqnarray}
\log \Big(1+\sum_{i \in {\cal I}, i\not = i^*}\exp \sum_{D\in {\cal D}^*_i}\theta_D(i_D)\Big)&=&
\sum_{l=1}^k\log \Big(1+\sum_{i \in {\cal I}_{P_l}, i\not = i^*}\exp \sum_{D\in {\cal D}^*_i}\theta^{P_l}_D(i_D)\Big)\nonumber\\
&&\hspace{2cm}
-\sum_{l=2}^k\log \Big(1+\sum_{i \in {\cal I}_{S_l}, i\not = i^*}\exp \sum_{D\in {\cal D}^*_i}\theta^{S_l}_D(i_D)\Big)\nonumber
\end{eqnarray}
it follows immediately that (\ref{pithetap}) is verified.

We note that as the restriction of $s$ to ${\cal D}^{P_l}, i_D\in {\cal I}^*_D$, the coefficients of $(\theta^{P_l}(i_D), D\in {\cal D}^{P_l}, i_D\in {\cal I}^*_D)$ in (\ref{pithetapl}) are  consistent across prime components and separators. The factorization of the induced conjugate prior on $p$ can be proved in a similar fashion.
\end{pff}

\vspace{2mm}

For given data $y$ with total count $n$, the posterior distribution of $\theta$ given $y$ will be
\begin{eqnarray}
\label{pithetapposterior}
\pi_G(\theta|s+y,\alpha+n)&=&\frac{\prod_{l=1}^k \pi_{P_l}(\theta^{P_l}|s^{P_l}+y^{S_l},\alpha+n)
}{\prod_{l=2}^k\pi_{S_l}(\theta^{S_l}|s^{S_l}+y^{S_l},\alpha+n)}
\end{eqnarray}
When comparing two models $G$ an d$G'$ the Bayes factor is the ratio of quantities of the type
$$\frac{I_{G'}(s,\alpha)}{I_{G'}(s,\alpha)}\;.$$
In the restricted class of decomposable graphical models, it is well-known that one can go from one decomposable graph to another through a succession of graphs that differ by only one edge. The additional edge can only belong to one clique in the new graph and as a consequence the Bayes factor affects only the graph induced by two cliques (see (37) in Dawid  and Lauritzen, 1993).
We are not aware of any such rule in the case of nondecomposable models. However, it is clear that the Bayes ratio will only involve the computation of the normalising constants for the subgraph induced by the prime components $P_l$ affected by the additional edge.

\section{Conclusion}
In this paper we have studied the conjugate prior for the log-linear parameters of discrete hierarchical log-linear models and its induced prior on the cell probability parameter $p$ thus extending the hyper Dirichlet which was the only form of the conjugate prior identified so far.

This prior has all the properties that one usually requires. As we have shown it, it has a moderate number 
of hyper-parameters precisely as many as there are log-linear parameters plus one. These hyperparameters are consistent across models. It is not difficult to translate prior knowledge into constraints for the hyper-parameters and constraints both in terms of the log-linear parameters and cell probabilities are consistent with prior beliefs, as illustrated in \S 5.1.

This prior has the additional property of being strong hyper Markov, thus leading to local updates for the computation of Bayes factors and it is also, of course, mathematically convenient since the prior and the posterior have the same form as the likelihood.
The conjugate  prior  should  therefore be one of the priors used for the study of contingency tables with a multinomial distribution for the cell counts. Though we have not mentioned it above, the translation of our results to the case of Poisson sampling is immediate.

\section{References}
\noindent Agresti, A. (1990). Categorical Data Analysis, {\it Wiley}.
\vspace{2mm}

\noindent Clyde, M. and George, E.I. (2004), Model uncertainty, {\it Statistical Science}, {\bf 19}, 81-94.
\vspace{.2cm}

\noindent Consonni, G. and Leucari, V. (2006). Reference priors for discrete graphical models, {\it Bka}, {\bf 93}, 23-40.
\vspace{.0cm}

\noindent  Darroch, J.N., Lauritzen, S.L. and Speed, T.P. (1980). Markov Fields and log-linear models for contingency tables, {\it Ann. Statist.}, {\bf 8}, 522-539.
\vspace{.2cm}

\noindent Dellaportas, P. and Forster, J.J. (1999). Markov chain Monte Carlo model determination for hierarchical and graphical log-linear models, {\it Bka}, {\bf 86}, 615-633.
\vspace{2mm}

\noindent Diaconis, P. \& Ylvisaker, D. (1979). Conjugate priors for exponential families. {\it Ann. Statist.}, {\bf 7}, 269-81.
\vspace{.2cm}

\noindent Hook, Albright and Cross (1980). Use of Bernoulli census and log-linear methods for estimating the prevalence of spina bifida in live births and the completeness of vital records reports in New York State. {\it Amer. J. Epidemiology}, {\bf 112}, 750-758.
\vspace{.2cm}

\noindent King, R. and Brooks, S.P., (2001). Prior induction for log-linear models for  general contingency table analysis, {\it Ann. Stat.}, {\bf 29}, 715-747.
\vspace{.2cm}  

\noindent Lauritzen, S. L. (1996).  {\sl Graphical Models.} Oxford University Press.
\vspace{.2cm}

\noindent Madigan, D. and Raftery, A.E. (1994). Model selection and accounting for model uncertainty in graphical models using Occam's window. {\it J.A.S.A.}, {\bf 89}, 1535-1546.
\vspace{.2cm}

\noindent Perks, W. (1947). Some observations on inverse probability including a new indifference rule. {\it J. Insti. Actuar.}. {\bf 73}, 285-334.
\vspace{.2cm}

\section{Appendix}
\subsection{Proof of Lemma \ref{lemmajacobian}}
We will first give the proof in the case where $G$ is the simple decomposable graph
$a-- b-- c$ and the data is bivariate. We will then sketch the proof for the general case of an arbitrary decomposable graph and discrete data. 

In the particular case of bivariate data, $p^G$ in Lemma \ref{lemmajacobian} becomes
\begin{equation}
\label{chosenp}
p^G=(p_D^{C_i}, D\in {\cal D}_{C_i}\setminus (\cup_{j=2}^k {\cal D}_{S_j}),\;i=1,\ldots,k,\;\; p_D^{S_i}, D\in {\cal D}_{S_i},\;i=2,\ldots,k)
\end{equation}
and the Jacobian of the change of variables from $\theta=(\theta_D, D\in {\cal D})$ as given in (\ref{thetae})  to $p^G$ as given in (\ref{chosenp}) is
\begin{equation}
\label{hyperJ}
\Big|\frac{d \theta}{d p^G}\Big|^{-1}=\frac{\prod_{i=1}^k\prod_{ D\in {\cal D}^{C_i}_0}p_D^{C_i}}
{\prod_{i=2}^k\prod_{D\in {\cal D}^{S_i}_0}p_D^{S_i}}
\end{equation}

\noindent We are therefore going to first prove (\ref{hyperJ}) for the two-chain graph above.
In this case we have $C_1=\{a,b\}, C_2=\{b,c\}, S=\{b\},\;\;p^G=(p_a^{C_1},p_{ab}^{C_1}, p_c^{C_2},p_{bc}^{C_2},p_b^S)$ and
\begin{eqnarray}
&&e^{\theta_a}=\frac{p_{a}}{p_{\emptyset}},\;e^{\theta_b}=\frac{p_{b}}{p_{\emptyset}},\;e^{\theta_c}=\frac{p_{c}}{p_{\emptyset}},\;e^{\theta_a+\theta_{ab}}=\frac{p_{ab}}{p_{b}},\;e^{\theta_{bc}+\theta_c}=\frac{p_{bc}}{p_{b}},\nonumber
\end{eqnarray}
Moreover, since the multinomial distribution is Markov with respect to the graph $G$, we have
$$p_{abc}=\frac{p_{ab}p_{bc}}{p_b}\;\;\;\mbox{and}\;\;\;p_{ac}=\frac{p_ap_c}{p_{\emptyset}}.$$
Therefore
\begin{eqnarray}
&&\frac{p_{a}^{C_1}}{p_{\emptyset}^{C_1}}
=\frac{p_{a}+p_{ac}}{p_{\emptyset}+p_c}=\frac{p_{a}}{p_{\emptyset}}
\frac{(1+\frac{p_c}{p_{\emptyset}})}{(1+\frac{p_c}{p_{\emptyset}})}
=\frac{p_{a}}{p_{\emptyset}}
\label{v1}\\
&&\frac{p_{c}^{C_2}}{p_{\emptyset}^{C_2}}=\frac{p_{c}}{p_{\emptyset}}\label{v2}\\
&&\frac{p_{ab}^{C_1}}{p_{b}^{C_1}}=\frac{p_{ab}+p_{abc}}{p_b+p_{bc}}
=\frac{p_{ab}}{p_b}\frac{(1+\frac{p_{bc}}{p_b})}{(1+\frac{p_{bc}}{p_b})}
=\frac{p_{ab}}{p_{b}}\label{v3}\\
&&\frac{p_{bc}^{C_2}}{p_{b}^{C_2}}=\frac{p_{bc}+p_{abc}}{p_b+p_{ab}}
=\frac{p_{bc}}{p_b}\frac{(1+\frac{p_{ab}}{p_b})}{(1+\frac{p_{ab}}{p_b})}
=\frac{p_{bc}}{p_{b}}\label{v4}\\
&&\frac{p_{b}^S}{p_{\emptyset}^S}=\frac{p_b+p_{ab}+p_{bc}+\frac{p_{ab}p_{bc}}{p_b}}{p_{\emptyset}+p_a+p_c+\frac{p_ap_c}{p_{\emptyset}}}=\frac{p_{b}}{p_{\emptyset}}\frac{
(1+\frac{p_{ab}}{p_b})(1+\frac{p_{bc}}{p_b})}
{(1+\frac{p_a}{p_{\emptyset}})(1+\frac{p_c}{p_{\emptyset}})
}\label{v5}
\end{eqnarray}
We introduce the intermediate variables 
$$v_a=\frac{p_{a}^{C_1}}{p_{\emptyset}^{C_1}},\;\;v_b=\frac{p_{b}^S}{p_{\emptyset}^S},\;\;v_c=\frac{p_{c}^{C_2}}{p_{\emptyset}^{C_2}},\;\;v_{ab}=\frac{p_{ab}^{C_1}}{p_{b}^{C_1}},\;\;v_{bc}=\frac{p_{bc}^{C_2}}{p_{b}^{C_2}}.$$ 

From (\ref{v1}) to (\ref{v5}), we have
\begin{eqnarray}
&&v_a=e^{\theta_a},\;v_c=e^{\theta_c},\;v_{ab}=e^{\theta_a+\theta_{ab}},\;v_{bc}=e^{\theta_{bc}+\theta_c}\;\;\mbox{and}\;\;e^{\theta_b}=v_b\frac{(1+v_a)(1+v_c)}{(1+v_{ab})(1+v_{bc})}\nonumber
\end{eqnarray}
It is then immediate to see that
\begin{equation}
\Big|\frac{d \theta}{d v}\Big|=\prod_{D\in {\cal D}}v_D^{-1}=\frac{p_{\emptyset}^{C_1}p_{\emptyset}^{S}p_{\emptyset}^{C_2}p_{b}^{C_1}p_{b}^{C_2}
}
{p_{a}^{C_1}p_{ab}^{C_1}p_{b}^{S}p_{c}^{C_2}p_{bc}^{C_2}
}
\end{equation}
Moreover, since $p_b^S=p_b^{C_1}+p_{ab}^{C_1}=p_b^{C_2}+p_{bc}^{C_2}$, then
$p_{\emptyset}^{C_1}=1-p_a^{C_1}-p_b^{C_1}-p_{ab}^{C_1}=1-p_a^{C_1}-p_b^S$ and similarly $p_{\emptyset}^{C_2}=1-p_c^{C_2}-p_b^S$. Therefore
\begin{equation}
v_a=\frac{p_a^{C_1}}{1-p_a^{C_1}-p_b^S},\;v_c=\frac{p_c^{C_2}}{1-p_c^{C_2}-p_b^S},\;v_{ab}=\frac{p_{ab}^{C_1}}{p_b^S-p_{ab}^{C_1}},\;v_{bc}=\frac{p_{bc}^{C_2}}{p_b^S-p_{bc}^{C_2}},\; v_b=\frac{p_b^S}{1-p_b^S}
\end{equation}
and the matrix of the Jacobian $\Big|\frac{dv}{dp^G}\Big|$ is
\begin{equation}
\frac{dv}{dp^G}=\pmatrix{\frac{1-p_b^S}{(1-p_a^{C_1}-p_b^S)^2}&0&0&0&0\cr
0&\frac{p_b^S}{(p_b^{C_1})^2}&0&0&0\cr
*&*&\frac{1}{(1-p_b^S)^2}&*&*\cr
0&0&0&\frac{p_b^S}{(p_b^{C_2})^2}&0\cr
0&0&0&0&\frac{1-p_b^S}{(1-p_c^{C_2}-p_b^S)^2}}
\end{equation}
The Jacobian is equal to the product of the diagonal elements and since $1-p_b^S=p_{\emptyset}^S$
\begin{equation}
\Big|\frac{dv}{dp^G}\Big|=
\frac{(p_{\emptyset}^S)^2(p_b^S)^2}
{(p_{\emptyset}^{C_1})^2(p_b^{C_1})^2(p_{\emptyset}^S)^2(p_b^{C_2})^2(p_{\emptyset}^{C_2})^2
}=\frac{(p_b^S)^2}
{(p_{\emptyset}^{C_1})^2(p_b^{C_1})^2(p_b^{C_2})^2(p_{\emptyset}^{C_2})^2
}
\end{equation}
Therefore
\begin{equation}
\Big|\frac{d\theta}{dp_G}\Big|=\frac{p_{\emptyset}^Sp_b^S
}
{p_{\emptyset}^{C_1}p_{a}^{C_1}p_{ab}^{C_1}p_{b}^{C_1}
p_{\emptyset}^{C_2}p_{b}^{C_2}p_{bc}^{C_2}p_{c}^{C_2}
}\;,
\end{equation}
which proves the lemma for the simple two-link  chain graph considered.

\noindent For a general decomposable graph with bivariate data, if we write $S=\cup_{i=2}^kS_i$, then, the intermediate variables will be
$$v=(\frac{p_D^{C_i}}{p_{D\cap S}^{C_i} }, D\in {\cal D}_{C_i}\setminus (\cup_{j=2}^k {\cal D}_{S_j}),\;i=1,\ldots,k,\;\; \frac{p_D^{S_i}}{p_{\emptyset}^{S_i}}, D\in {\cal D}_{S_i},\;i=2,\ldots,k)\;$$
and the proof will follow the same lines as above.

\noindent In the case of discrete data, the proof follows the same line as the proof above with the following substitutions. For $D\in {\cal D}$,
\begin{eqnarray}
&&\theta_D\;\;\mbox{becomes}\;\;(\theta_D(i_D), i_D\in {\cal I}_D^*)\nonumber\\
&&p_D;\;\mbox{becomes}\;\;(p(i_D), i_D\in {\cal I}_D^*)\nonumber\\
&&p_D^{C_l}\;\;\mbox{and}\;\;p_D^{S_l}\;\;\mbox{become} \;\;p^{C_i}(i_D)\;\;\mbox{and}\;\;p^{S_l}(i_D)\;\;\mbox{respectively},\;\; i_D\in {\cal I}_D^*\nonumber\\
&&v_D\;\;\mbox{becomes}\;\;(v_D(i_D), i_D\in {\cal I}_D^*)\nonumber\\
&&p_{\emptyset}^{C_l}=1-\sum_{D\in {\cal D}}p_D^{C_l};\;\mbox{becomes}\;\;p_{\emptyset}^{C_l}=1-\sum_{D\in {\cal D}}\sum_{i_D\in {\cal I}_D^*}p^{C_l}(i_D)\;\;\mbox{and similarly for}\;\; p_D^{S_l}\nonumber
\end{eqnarray}

\subsection{Proof of Lemma \ref{F}}
For ease of notation, we will give the proof of the lemma in the case of binary data. Since for each $C\in {\cal D}$ and $H\in {\cal E}$ there is only one cell in ${\cal I}^*_C$ and ${\cal I}_H^*$ respectively, we will adopt the notation 
$$F_{C,H}=F(i_C,j_H),\;\;C\in {\cal D}, H\in {\cal E}\;.$$
Let us first prove that if $H$ is decomposable, then (\ref{sum}) is true.
We proceed by induction on the number $k$ of cliques of $H$. Let ${\cal C}=\{C_1,\ldots,C_k\}$ be a perfect ordering of the cliques of $H$.

\noindent If $H$ is complete, that is $k=1$, we consider two cases, the case where $|H|$ is even and the case where it is odd.  For $|H|=2p, p\in {\bf N}$, there are $n_e=\sum_{k=1}^{p} \pmatrix{|H|\cr 2k}$ nonempty subsets of $H$ of even cardinality  and $n_o=\sum_{k=0}^{p-1} \pmatrix{|H|\cr 2k+1}$ subsets of odd cardinality. 
Therefore 
$$\sum_{C\subseteq_G H}(-1)^{|C|-1}=\sum_{k=1}^{2p}\pmatrix{2p\cr k}(-1)^{k+1}=(1-1)^{2p}-\pmatrix{2p\cr 0}(-1)^1=1$$
 and (\ref{sum}) is verified. We omit the proof for the case $|H|=2p-1$ which is parallel to that of the previous case. Therefore (\ref{sum}) is verified for $k=1$.

\noindent Let us now assume that $H$ is decomposable but not complete, that is $k>1$ and let us assume that (\ref{sum}) is true for any decomposable subset with $k-1$ cliques. It is well-know from the theory of decomposable graphs that, if we write $H_{k-1}=\cup_{j=1}^{k-1}C_j$, then $H=H_{k-1}\cup (C_k\setminus S_k)$ where $S_k=H_{k-1}\cap C_k$ is the $k$-th minimal separator in $H$. Therefore we have
\begin{equation}
\label{cminuss}
\sum_{C\subseteq_G H}F_{C,H}=\sum_{C\subseteq_G H_{k-1}}F_{C,H}+\Big(\sum_{C\subseteq_G C_k}F_{C,H}-\sum_{C\subseteq_G S_{k}}F_{C,H}\Big)\;.
\end{equation}
The first term on the right hand side of (\ref{cminuss}) is equal to 1 by our induction assumption while each one of the two other terms is also equal to 1 because both $C_k$ and $S_k$ are complete and therefore (\ref{sum}) is also verified for decomposable $H$.

\noindent Let us now prove that if $H$ is not decomposable and connected, $\sum_{C\subseteq_G H}F_{C,H}$ cannot be equal to 1.  If $H$ is not connected and its connected components $H^{(1)},\ldots, H^{(l)}$,  for some $l\geq 2$, are all decomposable, we clearly have $$\sum_{C\subseteq_G H}F_{C,H}=\sum_{j=1}^l\Big(\sum_{C\subseteq_G H^{(j)}}F_{C,H^{(j)}}\Big)\not =1.$$ 
 If $H$ is not connected and  its components are not all decomposable, this  implies that there is a nondecomposable subset $F_1$ of $G$ which can be separated from another subset $F_2$ of $G$ but this contradicts our assumption that $G$ is a prime component of $G$. So, this case does not occur.

\noindent If $H$ is not decomposable and connected, consider its set of cliques $\{C_1,\ldots,C_k\}$. Since $H$ is not decomposable, there is no perfect ordering of the cliques and therefore for any given ordering, there exist a nonempty subset $ {\cal Q}\subseteq \{3,\ldots,k\}$ such that for $j\in {\cal Q}$, there is no 
$i< j$ in the given ordering of the cliques of $H$ with 
$S_{j}=C_{j}\cap (\cup_{l=1}^{j-1}C_l)\subseteq C_i$ and therefore
$$S_{j}=C_{j}\cap (\cup_{l=1}^{j-1}C_l)=\oplus_{l=1}^{s_j}S_{j_{l}},\;\;\; 2\leq s_j\leq j-1$$
where the $S_{j_{l}}$ can be chosen to be disjoints, with $S_{j_{l}}\subseteq C_{j}\cap C_m$ for some $ m\in \{1,\ldots, j-1\}$.

For $j\in \overline{\cal Q}=\{2,\ldots, k\}\setminus {\cal Q}$, there exists $i<j$ in the given ordering of the cliques of $H$ such that $S_j\subseteq C_i$. Therefore
\begin{eqnarray}
\sum_{C\subseteq_G H}F_{C,H}&=&\sum_{C\subseteq_G C_1}F_{C,H}
+\sum_{j\in \overline{\cal Q}}\Big(\sum_{C\subseteq_G C_{j}}F_{C,H}-\sum_{C\subseteq_G S_{j}}F_{C,H}\Big)\label{l1}\\
&&+\sum_{j\in {\cal Q}}\Big(
\sum_{C\subseteq_G C_{j}}F_{C,H}-\sum_{l=1}^{s_j}\sum_{C\subseteq_G S_{j_l}}F_{C,H}\Big)\label{lq}\;.
\end{eqnarray}

The sums $\sum_{C\subseteq_G U}F_{C,H},\; U=C_1,C_j, S_j, j\in\overline{\cal Q}$ are all equal to 1 since each of $C_1, C_j, S_j, j\in\overline{\cal Q}$ are complete and connected and therefore the right hand side of 
 (\ref{l1}) is equal to 1. For the same reason, on line (\ref{lq}),
for $U=C_j, S_{j_l}, j\in {\cal Q}, l=1,\ldots,s_j, \sum_{C\subseteq_G U}F_{C,H}=1$. Since $s_j\geq 2,$
$$ \sum_{C\subseteq_G C_{j}}F_{C,H}-\sum_{l=1}^{s_j}\sum_{C\subseteq_G S_{j_l}}F_{C,H}\leq -1,\;\;j\in {\cal Q}$$ 
and therefore the sum on line (\ref{lq}) is less than or equal to $-|{\cal Q}|$. It follows that 
$$\sum_{C\subseteq_G H}F_{C,H}\leq 0$$ and in particular it cannot be equal to 1. The lemma is now proved.

\subsection{ Proof of Lemma \ref{jacob}}
Here again, we will give the proof of the lemma for binary data and we will use the notation of \S 2.3. To shorten notation, we will write $E\subseteq_G F$ to indicate that $E\subseteq F$ and $E\in {\cal D}$.

\noindent It is more convenient to compute $|\frac{d \theta}{d p}|$, express it in function of $\theta$ and take its inverse. From the expression (\ref{pe}) of $p_D, D\in {\cal D}$, we have 
\begin{eqnarray}
\frac{d p_D}{d \theta_D}&=&\frac{e^{\sum_{E\subseteq_G D}\theta_E}
}{
1+\sum_{F\in {\cal E}}e^{\sum_{E\subseteq_G F}\theta_E}}-
\frac{(e^{\sum_{E\subseteq_G D}\theta_E})\sum_{ F\supseteq D, F\in {\cal E}}e^{\sum_{E\subseteq_G F}\theta_E}
}
{(1+\sum_{F\in {\cal E}}e^{\sum_{D\subseteq_G F}\theta_D})^2}\nonumber\\
&=&\frac{e^{\sum_{E\subseteq_G D}\theta_E}
}{
1+\sum_{F\in {\cal E}}e^{\sum_{E\subseteq_G F}\theta_E}}
\Big(1-\frac{\sum_{F\in {\cal E}, F\supseteq D}e^{\sum_{E\subseteq_G F}\theta_E}
}
{1+\sum_{F\in {\cal E}}e^{\sum_{D\subseteq_G F}\theta_D}}\Big)\nonumber\\
&=&p_D(1-\sum_{F\in {\cal E}, F\supseteq D}p_F)\label{dpdd}\\
\frac{d p_D}{d \theta_C}&=&-
\frac{(e^{\sum_{E\subseteq_G D}\theta_E})\sum_{F\in {\cal E}, F\supseteq C}e^{\sum_{E\subseteq_G F}\theta_E}
}
{(1+\sum_{F\in {\cal E}}e^{\sum_{D\subseteq_G F,}\theta_D})^2},\;\;\;\;\;\;\; C\in {\cal D}, C\not \subset D\nonumber\\
&=&-p_D\sum_{F\in {\cal E}, F\supseteq C}p_F\;.\label{dpdc}\\
\frac{d p_D}{d \theta_C}&=&\frac{e^{\sum_{E\subseteq_G D}\theta_E}
}{
1+\sum_{F\in {\cal E}}e^{\sum_{E\subseteq_G F}\theta_E}}-
\frac{(e^{\sum_{E\subseteq_G D}\theta_E})\sum_{F\in {\cal E}, F\supseteq C}e^{\sum_{E\subseteq_G F}\theta_E}
}
{(1+\sum_{F\in {\cal E}}e^{\sum_{D\subseteq_G F}\theta_D})^2} ,\;\;\;\;C\subset D, C\not = D,\nonumber\\
&=&p_D(1-\sum_{F\in {\cal E}, F\supseteq C}p_F)\label{dpdcn}
\end{eqnarray}

We fix an arbitrary order of the elements of ${\cal D}$. From (\ref{dpdd}), (\ref{dpdc}) and (\ref{dpdcn}), it follows that 
the matrix of the Jacobian is such that the column of partial derivatives of $p_D$ is the vector with $C$-component
$$p_D\Big(\delta_{\{F\subseteq D\}}(C)-\sum_{F\in {\cal E}, F\supseteq C}p_F\Big),\;\;C\in {\cal D}\;,$$
where $\delta_{\{F\subseteq D\}}(C)$ is equal to 1 if $C\subseteq D$ and is equal to 0 otherwise. We note first that $p_D$ is common to all components of the $D$ column and therefore 
\begin{equation}
\label{product}
J=\det A\prod_{D\in {\cal D}}p_D
\end{equation}
where $A$ is the $|{\cal D}|\times|{\cal D}|$ matrix with entries
$$A_{C,D}=\delta_{\{F\subseteq D\}}(C)-\sum_{F\in {\cal E}, F\supseteq C}p_F,\;\;\; C, D \in {\cal D}.$$
We note next that in the rows $r_C$ corresponding to $C\in {\cal D}$ maximal with respect to inclusion, the $C$ entry, on the matrix diagonal, is the only entry such that 
$$\delta_{\{F\subseteq D\}}(C)\not = 0\;,$$
and therefore, for $C$ maximal, we can write
\begin{equation}
\label{cmax}
r_C=e_C-(\sum_{F\in {\cal E}, F\supseteq C}p_F){\bf 1}^t\;,
\end{equation}
where $e_C$ is the ${\cal D}$-dimensional row vector with components all equal to 0 except for the $C$ component, and ${\bf 1}^t$ is the ${\cal D}$-dimensional row vector with all its components equal to 1. 

\noindent We finally note that if $C_1\subset C_2$ for $C_1,C_2$ in ${\cal D}$, then 
$$\{F\in {\cal E}, F\supseteq C_2\}\subset \{F\in {\cal E}, F\supseteq C_1\}$$
and therefore if, in the matrix $A$, for $C\in {\cal D}$  not maximal with respect to inclusion, we replace the row $r_C$ by
$$\tilde{r}_C=r_C+\sum_{F\supset C, F\in {\cal D}}(-1)^{|F\setminus C|}r_F$$
we have
\begin{equation}
\label{cnotmax}
\tilde{r}_C=e_C-\Big(\sum_{F\supseteq C, F\in {\cal D}}(-1)^{|F\setminus C|}(\sum_{H\in {\cal E}, H\supseteq F}p_H)\Big){\bf 1}^t\;.
\end{equation}
The determinant of $A$ is clearly equal to the determinant of  the matrix $\tilde{A}$ obtained from $A$ by replacing $r_C$ by $\tilde{r_C}$ whenever $C\in {\cal D}$ is not maximal with respect to inclusion. Using (\ref{cmax}) and (\ref{cnotmax}), we have
$$\tilde{A}=I_{|{\cal D}|}-U{\bf 1}^t\;$$
where $U$ is the column vector $U=(\sum_{F\supseteq C, F\in {\cal D}}(-1)^{|F\setminus C|}(\sum_{H\in {\cal E}, H\supseteq F}p_H),\; C\in {\cal D}).$ It is well-known that 
$$\det \tilde{A}=1-{\bf 1}^tU$$
Therefore
\begin{eqnarray}
\det \tilde{A}&=&1-\sum_{C\in {\cal D}}\Big(\sum_{D\supseteq C, D\in {\cal D}}(-1)^{|D\setminus C|}(\sum_{H\in {\cal E}, H\supseteq D}p_H)\Big)\nonumber\\
&=&1-\sum_{H\in {\cal E}}p_H(\sum_{\{D\subseteq H, D\in {\cal D}\}}\sum_{\{C\subseteq D, C\in {\cal D}\}}(-1)^{|D\setminus C|})\label{stop}
\end{eqnarray}
According to (\ref{fich}), the coefficients of $p_H$ in the expression above are the sum of the entries  $F_{D,H}=\sum_{\{C\subseteq_G D\}}(-1)^{|D\setminus C|})$ in the column $H$ of $F$. Moreover, by Lemma \ref{F}, this sum $\sum_{D\subset H, D\in {\cal D}}F_{D,H}$ is  equal to 1 if and only if  $H\in {\cal E}$ is decomposable, connected and nonempty. Since $1=\sum_{F\in {\cal E}_0}p_F=\sum_{F\in {\cal U}_0}p_F+\sum_{F\not \in {\cal U}_0}p_F$, we can write
\begin{eqnarray}
\det \tilde{A}&=&\sum_{F\in {\cal U}_0}p_F+\sum_{F\not \in {\cal U}_0}p_F-\sum_{H\in {\cal E}}p_H(\sum_{\{D\subseteq H, D\in {\cal D}\}}F_{D,H})\nonumber\\
&=&\sum_{F\in {\cal U}_0}p_F+\sum_{F\not \in {\cal U}_0}p_F-\sum_{H\not \in {\cal U}_0}p_H-\sum_{H\in {\cal U}_0}p_H(\sum_{\{D\subseteq H, D\in {\cal D}\}}F_{D,H})\nonumber\\
&=&-\sum_{H\in {\cal U}_0}p_H\Big((\sum_{\{D\subseteq H, D\in {\cal D}\}}F_{D,H})-1\Big)\nonumber\\
&=&-\sum_{H\in {\cal U}_0}a_Hp_H \label{productp}
\end{eqnarray}
From (\ref{product}) and (\ref{productp}), we derive the first expression for $J$ in (\ref{jacobianptheta}). The other expressions are deduced by replacing the different $p_F$ by their expression with respect to $(\theta_D,D\in {\cal D})$.
\begin{remark}
\label{notgraphical}
When the model is not specified to be graphical but is more generally hierarchical, the proof above holds up to (\ref{stop}).
\end{remark}
\end{document}